\documentclass[a4paper,12pt]{article}
\usepackage{latexsym,amsfonts}

\title{Long-time dynamics for the 3-D viscous primitive equations of large-scale moist atmosphere}

\author{Boling Guo
and Daiwen Huang\thanks{Corresponding author's E-mail:
hdw55@tom.com.}
\\ \small  Institute of Applied Physics and Computational
Mathematics, Beijing 100088, China\\}
\date{}

\begin{document}

\maketitle

\begin{minipage}{12cm}
{\small {\bf Abstract.}\quad In this paper, we consider the
initial boundary value problem for the three-dimensional viscous
primitive equations of large-scale moist atmosphere which are used
to describe the turbulent behavior of long-term weather prediction
and climate changes. First, we  obtain the existence and
uniqueness of global strong solutions of the problem. Second, by
studying the long-time behavior of  strong solutions, we construct
a (weak) universal attractor $\mathcal{A}$ which captures all the
trajectories.
\smallskip


\smallskip

}
\end{minipage}

\section{Introduction}
In order to understand the mechanism of long-term weather
prediction and climate changes, one can study the mathematical
equations and models governing the motion of the atmosphere as the
atmosphere is a specific compressible fluid (see, e.g., \cite{H,
P}). V. Bjerkness, one of the pioneers of meteorology, said that
the weather forecasting can be considered as an initial boundary
value problem in mathematical physics. In 1922, Richardson
initially introduced the so-called primitive atmospheric equations
which consisted of the hydrodynamic, thermodynamic equations with
Coriolis force, cf. \cite{R}. At that time, the primitive
atmospheric equations were too complicated to be studied
theoretically or to be solved numerically. To overcome this
difficulty, some simple numerical models were introduced, such as
the barotropic model formulated by Neumann etc. in \cite{CFN} and
the quasi-geostrophic model introduced by Charney and Philips in
\cite{CP}. The 2-D and 3-D quasi-geostrophic models have been the
subject of analytical mathematical study, cf., e.g., \cite{B, C,
CMT1, CMT2, CW, EM, M, W2, W3, W4} and references therein.

Due to the considerable improvement in computer capacity and the
development of atmospheric science, some mathematicians began to
consider the primitive equations of atmosphere in the past two
decades (see, e.g., \cite{LTW1, LTW2, LTW3, TZ} and references
therein). In \cite{LTW1}, by introducing viscosity terms and using
some technical treatment, Lions, Temam and Wang obtained a new
formulation of the primitive equations of large-scale dry
atmosphere which was amenable to mathematical study. In a
p-coordinate system, the new formulation of the primitive
equations is a little similar to  Navier-Stokes equations of
incompressible fluid. By the methods used to solve Navier-Stokes
system in \cite{L}, they obtained the existence of global weak
solutions of the initial boundary value problem for the new
formulation of the primitive equations. Moreover, under the
assumptions that there exists a unique global strong solution for
the problem with vertical viscosity and that $H^1$-norm of the
strong solution is uniformly bounded about $t$, they established
some physically relevant estimates for the Hausdorff and fractal
dimensions of the attractor of the primitive equations with
vertical viscosity. Without those assumptions, by the Trajectory
Attractors Theory due to Vishik and Chepyzhov ( cf. \cite{CV,VC}),
we obtained the existence of trajectory attractors for the
large-scale moist atmospheric primitive equations in \cite{GH}. By
taking advantage of the geostrophic balance and other geophysical
consideration, several intermediate models have been the subject
of studying the long-time dynamics and global attractors in order
to describe the mechanism of long-term weather prediction and
climate dynamics (see, e. g., \cite{CT1, CT2, JC1, W1, W2} and
references therein).

In recent years, there were some mathematicians who considered the
existence of strong solutions for the three-dimensional viscous
primitive equations of large-scale atmosphere and ocean (see, e.
g., \cite{C1, CT1, CT2, CT3, GM, HTZ, TSW, TZ} and references
therein). In \cite{GM}, Guill\'en-Gonz\'alez etc. obtained the
global existence of strong solutions to the primitive equations of
large-scale ocean by assuming that the initial data are small
enough, and also proved the local existence of strong solutions to
the equation for all initial data. In \cite{TZ}, Temam and Ziane
considered the local existence of strong solutions for the
primitive equations of the atmosphere, the ocean and the coupled
atmosphere-ocean. The papers \cite{CT1, CT2, CT3} are devoted to
considering  the non-dimensional Boussinesq equations or modified
models (see, e.g., \cite{P, TSW}). In \cite{CT1}, Cao and Titi
considered global well-posedness and finite-dimensional global
attractor to a 3-D planetary geostrophic model. The paper
\cite{CT3} is devoted to studying the global well-posedness for
the three-dimensional viscous primitive equations of large-scale
ocean. In \cite{CT3}, Cao and Titi developed a beautiful approach,
by which they obtained the fact that $L^6$-norm of the fluctuation
$\tilde{v}$ of horizontal velocity is bounded uniformly about the
time $t$. The estimate about $L^6$-norm of the fluctuation
$\tilde{v}$ is a key proof in \cite{CT3}. On the basis of the
results of \cite{CT3}, we obtain the existence of (weak) universal
attractors for the 3-D viscous primitive equations of  the
large-scale ocean in \cite{GH1}.

In the present paper we are interested in considering the
existence, uniqueness and long-time behavior of global strong
solutions to the initial boundary value problem of the new
formulation of large-scale moist atmospheric primitive equations
(the problem is denoted by (IBVP) and will be given in the section
2). Our main results are Theorem 3.1, Theorem 3.2, Proposition 3.3
and Theorem 3.4. First, we obtain the global well-posedness of the
problem (IBVP). Second, by studying the long-time behavior of the
strong solution, we prove $H^1$-norm of the strong solution is
uniformly bounded about $t$, and  also the corresponding semigroup
$\{S(t)\}_{t\geq 0}$ possesses a bounded absorbing set $B_\rho$ in
$V$ (the definition of the space $V$ will be given in subsection
4.1), by which we construct a (weak) universal attractor
$\mathcal{A}$. Since the global well-posedness of the 3-D
incompressible Navier-Stokes system is still open, by Theorem 3.1
and Theorem 3.2, we prove rigorously in mathematics that the new
formulation of large-scale moist atmospheric primitive equations
is simpler than the incompressible Navier-Stokes system, which is
consistent with the physical point of view.

Inspired by the methods used in \cite{CT3}, we prove global
well-posedness of (IBVP). However, there are two differences
between the results of our paper and that in \cite{CT3}. On one
hand, the new formulation of the large-scale moist atmospheric
equations is more complicated than the model studied in
\cite{CT3}. If we let $a=0$, the model considered in this paper is
similar to that in \cite{CT3}. On the other hand, we have studied
the long-time dynamics and the existence of the (weak) universal
attractors for the large-scale moist atmospheric primitive
equations. Here, the results about the universal attractors in
this paper is more stronger than that in \cite{GH}.

In order to study the long-time behavior of strong solutions, we
must make three key estimates. First, we must make estimates about
$L^3$-norm of the temperature $T$ and the fluctuation $\tilde{v}$
of horizontal velocity $v$ before we study the long-time behavior
of strong solutions by the Uniform Gronwall Lemma, without which
we only obtain the global well-posedness of (IBVP). Second, we
ought to make estimates about $L^4$-norm of $\tilde{v}$, $T$ and
the mixing ratio of water vapor in the air $q$. If we only made
estimates about $L^6$-norm of $\tilde{v},\ q,\ T$ as that in
\cite{CT3}, we could not study long-time behavior of strong
solutions and could not obtain a stronger result than the
uniqueness of strong solutions to (IBVP). Third, since the moist
atmospheric equations are more complicated than the oceanic
primitive equations, we have to make estimates about $\partial_\xi
T,\
\partial_\xi q$ before we prove $H^1$-norm of $v,\ T,\ q$ is bounded.

The paper is organized as follows: In section 2, we  pose the
primitive equations of large-scale moist atmosphere. Main results
of this paper are formulated in section 3. In section 4, we give
our working spaces and some preliminaries. We prove main results
of our paper in sections 5, 6, 7.

\section{The three-dimensional viscous primitive equations of large-scale moist atmosphere}
The three-dimensional viscous primitive equations of large-scale
moist atmosphere in the pressure coordinate system(for details, we
refer the reader to \cite{GH, JC2, LTW1, LTW2} and references
therein) is written as
$$\frac{\partial v}{\partial t}+\nabla_v v+\omega\frac{\partial v}{\partial \xi}+
\frac{f}{R_0}k\times v+\hbox{grad}\Phi-\frac{1}{Re_1}\triangle
v-\frac{1}{Re_2}\frac{\partial^2 v}{\partial \xi^2}=0,
\eqno{(2.1)}$$
$$\hbox{div}v+\frac{\partial \omega}{\partial\xi}=0,
\quad\quad\quad\quad\quad\quad\quad\quad\quad
\quad\quad\quad\quad\quad
\quad\quad\quad\quad\quad\quad\quad\eqno{(2.2)}$$
$$\frac{\partial \Phi}{\partial\xi}+\frac{bP}{p}(1+aq)T=0,
\quad\quad\quad\quad\quad\quad\quad\quad\quad
\quad\quad\quad\quad\quad\quad\quad\quad\ \eqno{(2.3)}$$
$$\frac{\partial T}{\partial t}+\nabla_v T+\omega\frac{\partial T}{\partial \xi}
-\frac{bP}{p}(1+aq)\omega-\frac{1}{Rt_1}\triangle
T-\frac{1}{Rt_2}\frac{\partial^2 T}{\partial \xi^2}=Q_1, \
\eqno{(2.4)}$$
$$\frac{\partial q}{\partial t}+\nabla_v q+\omega\frac{\partial q}{\partial \xi}
-\frac{1}{Rq_1}\triangle q-\frac{1}{Rq_2}\frac{\partial^2
q}{\partial \xi^2}=Q_2, \quad\quad\quad\quad\quad\quad\quad\quad\
\eqno{(2.5)}$$ where the unknown functions are $v,\ \omega,\
\Phi,\ q,\ T$, $v=(v_{\theta},v_{\varphi})$  the horizontal
velocity, $\omega$ vertical velocity in $p-$ coordinate system,
$\Phi$ the geopotential, $q$ the mixing ratio of water vapor in
the air, $T$ temperature, $f=2\cos \theta$ Coriolis parameter,
$R_0$ the Rossby number, $k$ vertical unit vector, $Re_1,\ Re_2$,
$Rt_1,\ Rt_2$, $Rq_1,\ Rq_2$ Reynolds numbers, $P$ an approximate
value of pressure at the surface of the earth,  $p_0$ pressure of
the upper atmosphere and $p_0>0$, the variable $\xi$ satisfying
$p=(P-p_0)\xi+p_0$ ($0<p_0\leq p\leq P$),  $Q_1,\ Q_2$ given
functions on $S^2\times (0,1)$ (here we don't consider the
condensation of water vapor), $a$ a positive constant($a\approx
0.618$), $b$ a positive constant. The definitions of $\nabla_v v,\
\triangle v,\ \triangle T,$ $\ \triangle q,\ \nabla_v q,\ \nabla_v
T,\ \hbox{div}v,$ $\ \hbox{grad}\Phi$ will be given in the
subsection 4.1. The equations $(2.1)-(2.5)$ are called the 3-D
viscous primitive equations of the large-scale moist atmosphere.

The space domain of the equations $(2.1)-(2.5)$ is
$$\Omega=S^2\times (0,1),$$ where $S^2$ is two-dimensional unit sphere. The boundary
value conditions are given by
$$\xi=1(p=P):\ \frac{\partial v}{\partial \xi}=0,\ \omega=0,\ \frac{\partial T}{\partial \xi}=
\alpha_s(T_s-T),\ \frac{\partial q}{\partial \xi}= \beta_s(q_s-q),
\eqno{(2.6)}$$
$$\xi=0(p=p_0):\ \ \frac{\partial v}{\partial \xi}=0,\ \omega=0,\ \frac{\partial T}{\partial \xi}=0,\
\frac{\partial q}{\partial \xi}= 0,
\quad\quad\quad\quad\quad\quad\quad\quad\eqno{(2.7)}$$ where
$\alpha_s,\ \beta_s$ are positive constants, $T_s$ the given
temperature on the surface of the earth, $q_s$ the given mixing
ratio of water vapor on the surface of the earth. For simplicity
and without loss generality we assume that $T_s=0$  and $q_s=0$.
If $T_s\neq 0$ and $q_s\neq 0$, one can homogenize the boundary
value conditions for $T,\ q$ (cf., e.g., \cite{GH}).

Integrating $(2.2)$ and using the boundary conditions $(2.6),\
(2.7)$, we have
$$\omega(t;\theta,\varphi,\xi)=W(v)(t;\theta,\varphi,\xi)=\int^1_{\xi}\hbox{div}
v(t;\theta,\varphi,\xi') \ d\xi',\eqno{(2.8)}$$
$$\int^1_0\hbox{div}
v\ d\xi=0.\eqno{(2.9)}$$ Suppose that $\Phi_s$ is a certain
unknown function at the isobaric surface $\xi=1$. Integrating
$(2.3)$, we obtain
$$\Phi(t;\theta, \varphi, \xi)=\Phi_s(t;\theta, \varphi)+\int^1_{\xi}\frac{bP}{p}
(1+aq)T \ d\xi'.\eqno{(2.10)}$$  Then the equations $(2.1)-(2.5)$
can be written as
$$\frac{\partial v}{\partial t}+\nabla_v v+W(v)\frac{\partial v}{\partial \xi}+
\frac{f}{R_0}k\times v+\hbox{grad}\Phi_s+\int^1_{\xi}\frac{bP}{p}
\hbox{grad}[(1+aq)T] \ d\xi'$$
$$-\frac{1}{Re_1}\triangle
v-\frac{1}{Re_2}\frac{\partial^2 v}{\partial \xi^2}=0,
\quad\quad\quad\quad\quad \quad\quad\eqno{(2.11)}$$
$$\frac{\partial T}{\partial t}+\nabla_v T+W(v)\frac{\partial T}{\partial \xi}
-\frac{bP}{p}(1+aq)W(v)-\frac{1}{Rt_1}\triangle
T-\frac{1}{Rt_2}\frac{\partial^2 T}{\partial \xi^2}=Q_1,
\eqno{(2.12)}$$
$$\frac{\partial q}{\partial t}+\nabla_v q+W(v)\frac{\partial q}{\partial \xi}
-\frac{1}{Rq_1}\triangle q-\frac{1}{Rq_2}\frac{\partial^2
q}{\partial \xi^2}=Q_2, \quad\quad\quad\quad\quad
\quad\quad\quad\eqno{(2.13)}$$
$$\int^1_0\hbox{div}
v\ d\xi=0,\quad\quad\quad\quad\quad\quad\quad\quad
\quad\quad\quad\quad\quad\quad\quad\quad \quad\quad\quad
\quad\quad\quad\eqno{(2.14)}$$ where the definitions of
$\hbox{grad}[(1+aq)T],\ \hbox{grad}\Phi_s$ will be given in the
subsection 4.1. The boundary value conditions of the equations
$(2.11)-(2.14)$ are given by
$$\xi=1:\ \frac{\partial v}{\partial \xi}=0,\ \frac{\partial T}{\partial \xi}=
-\alpha_s T,\ \frac{\partial q}{\partial \xi}= -\beta_s q,
\eqno{(2.15)}$$
$$\xi=0:\ \frac{\partial v}{\partial \xi}=0,\ \frac{\partial T}{\partial \xi}=0,\
\frac{\partial q}{\partial \xi}=
0;\quad\quad\quad\quad\eqno{(2.16)}$$ and the initial value
conditions can be given as
$$U|_{t=0}=(v|_{t=0},T|_{t=0},q|_{t=0})=U_0=(v_0,T_0,q_0).\eqno{(2.17)}$$ We call $(2.11)-(2.17)$ as the
initial boundary value problem of the new formulation of the 3-D
viscous primitive equations of large-scale moist atmosphere, which
is denoted by (IBVP).

Now we define the fluctuation $\tilde{v}$ of horizontal velocity
and find the equations satisfied by $\tilde{v}$ and $\bar{v}$. By
integrating the momentum equation(2.11) with respect to $\xi$ from
$0$ to $1$ and using the boundary value conditions (2.15) and
(2.16), we get
$$\frac{\partial \bar{v}}{\partial t}+\displaystyle\int_{0}^{1}(\nabla_v v
+W(v)\frac{\partial v}{\partial
\xi})d\xi+\frac{f}{R_0}k\times\bar{v}+\hbox{grad}
\Phi_s+\displaystyle\int_{0}^{1}\displaystyle\int_{\xi}^{1}\frac{bP}{p}\hbox{grad}[(1+aq)T]d\xi^{'}d\xi
$$ $$-\frac{1}{Re_1}\triangle\bar{v}=0\ \ \hbox{in}\ S^2, \eqno(2.18)$$ \noindent
where $\bar{v}=\displaystyle\int_{0}^{1}vd\xi.$

Denote the fluctuation of the horizontal velocity by
$$\tilde{v}=v-\bar{v}.$$
We notice that
$$\bar{\tilde{v}}=\displaystyle\int_{0}^{1}\tilde{v}d\xi=0,
\ \  \nabla \cdot \bar{v}=0. \eqno(2.19)$$ By integration by parts
and (2.19), we have
$$\displaystyle\int_{0}^{1} W(v)\frac{\partial v}{\partial \xi}d\xi
=\displaystyle\int_{0}^{1}v\hbox{div}vd\xi=\displaystyle\int_{0}^{1}\tilde{v}\hbox{div}\tilde{v}d\xi,
\eqno(2.20)$$
$$\displaystyle\int_{0}^{1}\nabla_v vd\xi
=\displaystyle\int_{0}^{1}\nabla_{\tilde{v}}
\tilde{v}d\xi+\nabla_{\bar{v}}\bar{v}. \eqno(2.21)$$ From (2.18),
(2.20) and (2.21), we obtain
$$\frac{\partial \bar{v}}{\partial t}+\nabla_{\bar{v}}\bar{v}
+\overline{\tilde{v}\hbox{div}\tilde{v}+\nabla_{\tilde{v}}
\tilde{v}} +\frac{f}{R_0}k\times\bar{v}+\hbox{grad}\Phi_s
+\displaystyle\int_{0}^{1}\int_{\xi}^{1}\frac{bP}{p}\hbox{grad}[(1+aq)T]d\xi^{'}d\xi$$
$$-\frac{1}{Re_1}\triangle \bar{v}=0\ \ \hbox{in}\ S^2. \eqno(2.22)$$
Subtracting (2.22) from (2.11), we know that the fluctuation
$\tilde{v}$ satisfies the following equation and boundary value
conditions
$$\frac{\partial \tilde{v}}{\partial t}+\nabla_{\tilde{v}} \tilde{v}
+(\displaystyle\int_{\xi}^{1}\hbox{div}\tilde{v}d\xi^{'})\frac{\partial
\tilde{v}}{\partial \xi} +\nabla_{\tilde{v}}
\bar{v}+\nabla_{\bar{v}}
\tilde{v}-\overline{(\tilde{v}\hbox{div}\tilde{v}+\nabla_{\tilde{v}}
\tilde{v})}+\frac{f}{R_0}k\times\tilde{v}$$
$$+\displaystyle\int_{\xi}^{1}\frac{bP}{p}\hbox{grad}[(1+aq)T]d\xi^{'}
-\displaystyle\int_{0}^{1}\int_{\xi}^{1}\frac{bP}{p}\hbox{grad}[(1+aq)T]d\xi^{'}d\xi
\quad\quad\quad\quad$$
$$
-\frac{1}{Re_1}\triangle \tilde{v}-\frac{1}{Re_2}\frac{\partial^2
\tilde{v}}{\partial \xi^2}=0\ \ \hbox{in}\ \Omega,
\quad\quad\quad\quad\quad\quad\quad\quad\quad\quad\quad\quad\quad\quad\eqno(2.23)$$
$$\xi=1: \frac{\partial \tilde{v}}{\partial \xi}=0,\quad\quad\quad\quad\quad\quad\quad\quad
\quad\quad\quad\quad\quad\quad\quad\quad\quad\quad\quad\quad
\eqno(2.24)$$
$$\xi=0: \frac{\partial \tilde{v}}{\partial \xi}=0. \quad\quad\quad\quad\quad\quad\quad\quad
\quad\quad\quad\quad\quad\quad\quad\quad\quad\quad\quad\quad\eqno(2.25)$$

\section{Statements of main results}
Now we formulate our main results in the present paper.

\smallskip
\noindent{\bf Theorem 3.1} \quad Let $Q_1,\ Q_2\in H^1(\Omega),\
U_0=(v_0,T_0,q_0)\in V$. Then for any $\mathcal{T}>0$ given, there
exists a strong solution $U$ of the system (2.11)-(2.17) on the
interval $[0,\mathcal{T}]$, where the definition of strong
solutions to the system (2.11)-(2.17) will be given in subsection
5.1.

\smallskip
\noindent{\bf Theorem 3.2} \quad Let $Q_1,\ Q_2\in H^1(\Omega),\
U_0=(v_0,T_0,q_0)\in V$. Then for any $\mathcal{T}>0$ given, the
strong solution $U$ of the system (2.11)-(2.17) on the interval
$[0,\mathcal{T}]$ is unique. Moreover, the strong solution $U$ is
dependent continuously on the initial data.

\smallskip
\noindent{\bf Proposition 3.3} \quad If $Q_1,\ Q_2\in
H^1(\Omega),\ U_0=(v_0,T_0,q_0)\in V$, Then the global strong
solution $U$ of the system (2.11)-(2.17) satisfies $U\in
L^{\infty}(0,\infty;V)$. Moreover, the corresponding semigroup
$\{S(t)\}_{t\geq 0}$ possesses a bounded absorbing set $B_\rho$ in
$V$, i.e., for every bounded set $B\subset V$, there exists
$t_0(B)>0$ big enough such that
$$S(t)B\subset B_\rho, \ \hbox{for any}\ t\geq t_0,$$ where $B_\rho=\{U; \|U\|\leq \rho\}$ and $\rho$ is a
positive constant dependent on $\|Q_1\|_1,\ \|Q_2\|_1$.

\smallskip
\noindent{\bf Theorem 3.4} \quad The system (2.11)-(2.16)
possesses a (weak) universal attractor
$\mathcal{A}=\displaystyle\cap_{s\geq 0}
\overline{\displaystyle\cup_{t\geq s}T(t)B_\rho}$ that captures
all the trajectories, where the closures are taken with respect to
$V$-weak topology. The (weak) universal attractor $\mathcal{A}$
has the following properties:

(i) (weak compact) \quad $\mathcal{A}$ is bounded and weakly
closed in $V$;

(ii) (invariant) \quad for every $t\geq 0$,
$S(t)\mathcal{A}=\mathcal{A}$;

(ii) (attracting) \quad for every bounded set $B$ in $V$, the sets
$S(t)B$ converge to $\mathcal{A}$ with respect to $V$-weak
topology as $t\rightarrow +\infty $, i.e.,
$$\displaystyle\lim_{t\rightarrow +\infty}d_V^w(S(t)B,\ \mathcal{A})=0,$$
where the distance $d_V^w$ is induced by the $V$-weak topology.

\smallskip
\noindent{\bf Remark 3.5} \quad The (weak) universal attractor
$\mathcal{A}$ has the following additional properties:

(i) \quad By Rellich-Kondrachov Compact Embedding Theorem (cf.,
e.g., \cite{A}), we know that for any $1\leq p<6$ the sets $S(t)B$
converge to $\mathcal{A}$ with respect to the $L^p(\Omega)\times
L^p(\Omega)\times L^p(\Omega)\times L^p(\Omega)$-norm;

(ii) The (weak) universal attractor $\mathcal{A}$ is unique and is
connected with respect to $V$-weak topology.

\smallskip
\noindent{\bf Remark 3.6} \quad In the forthcoming paper, we shall
prove the (weak) universal attractor $\mathcal{A}$ have finite
fractal and Hausdorff dimensions in $L^2(\Omega)\times
L^2(\Omega)\times L^2(\Omega)\times L^2(\Omega)$.

\smallskip
\noindent{\bf Remark 3.7} \quad In comparison to the 3-D
incompressible Navier-Stokes equations, the 3-D viscous primitive
equations of large-scale moist atmosphere have not the time
derivative term of the vertical velocity $\omega=W(v)$. Therefore,
we can not prove that the bounded absorbing set $B_\rho$ in $V$ is
bounded in $H^2(\Omega)\times H^2(\Omega)\times H^2(\Omega)\times
H^2(\Omega)$ as in the case of 3-D incompressible Navier-Stokes
equations (for the Navier-Stokes equations, if there exists a
bounded absorbing set $B_\rho$ in $H^1_0(\Omega)\times
H^1_0(\Omega)\times H^1_0(\Omega)$, then one can prove $B_\rho$ is
bounded in $H^2(\Omega)\times H^2(\Omega)\times H^2(\Omega)$),
i.e., we can not prove that the universal attractor $\mathcal{A}$
is compact in $V$. In this sense, the primitive equations of the
large-scale atmosphere are more complicated than 3-D
incompressible Navier-Stokes equations.

\section{Preliminaries}
\subsection{Some function spaces}
Let $e_{\theta},\ e_{\varphi},\ e_{\xi}$ be the unit vectors in
$\theta,\ \varphi$ and $\xi$ directions of the space domain
$\Omega$ respectively,
$$e_{\theta}=\frac{\partial}{\partial \theta},\ \ e_{\varphi}
=\frac{1}{\sin \theta}\frac{\partial}{\partial \varphi},\ \
e_{\xi}=\frac{\partial}{\partial \xi}.$$ The inner product and
norm on $T_{(\theta,\varphi,\xi)}\Omega$ (the tangent space of
$\Omega$ at the point $(\theta,\varphi,\xi)$) are given by
$$(X,Y)=X\cdot Y=X_1Y_1+X_2Y_2+X_3Y_3,\ \ |X|=(X,X)^{\frac{1}{2}}$$
for $$X=X_1e_{\theta}+X_2e_{\varphi}+X_3e_{\xi},\
Y=Y_1e_{\theta}+Y_2e_{\varphi}+Y_3e_{\xi}\in
T_{(\theta,\varphi,\xi)}\Omega.$$

$L^p(\Omega):=\{h;\ h: \Omega\rightarrow \mathbb{R},\int_\Omega
|h|^p<+\infty\}$ with the norm
$|h|_p=(\int_\Omega|h|^p)^{\frac{1}{p}}$, $1\leq p<\infty$.
$\int_\Omega \cdot d\Omega$ and $\int_{S^2} \cdot dS^2$ are
denoted by $\int_\Omega \cdot $ and $\int_{S^2} \cdot $
respectively. $L^2(T\Omega|TS^2)$ is the first two components of
$L^2$ vector fields on $\Omega$ with the norm
$|v|_2=(\int_\Omega(|v_{\theta}|^2+|v_{\varphi}|^2))^{\frac{1}{2}}$,
where $v=(v_{\theta},v_{\varphi}):\Omega\rightarrow TS^2$.
$C^{\infty}(S^2)$ is the function space for all smooth functions
from $S^2$ to $\mathbb{R}$. $C^{\infty}(\Omega)$ is the function
space for all smooth functions from $\Omega$ to $\mathbb{R}$.
$C^{\infty}(T\Omega|TS^2)$ is the first two components of smooth
vector fields on $\Omega$.  $H^m(\Omega)$ is the Sobolev space of
functions which are in $L^2$, together with all their covariant
derivatives with respect to $e_{\theta},\ e_{\varphi},\ e_{\xi}$
of order $\leq m$, with the norm
$$\|h\|_m=(\int_\Omega(\displaystyle\sum_{1\leq k\leq m}
\displaystyle\sum_{i_j=1,2,3;j=1,\cdot\cdot\cdot,k}|\nabla_{i_1}
\cdot\cdot\cdot\nabla_{i_k}h|^2+|h|^2))^{\frac{1}{2}},$$ where
$\nabla_1=\nabla_{e_{\theta}},\ \nabla_2=\nabla_{e_{\varphi}},\
\nabla_3=\nabla_{e_{\xi}}=\frac{\partial}{\partial \xi}$ (the
definitions of $\nabla_{e_{\theta}},\nabla_{e_{\varphi}}$ will be
given later). $H^m(T\Omega|TS^2)= \{v;\
v=(v_{\theta},v_{\varphi}): \Omega\rightarrow TS^2,\ \|v\|^m_m<
+\infty\},$ the norm of which is similar to that of $H^m(\Omega)$,
that is , in the above formula of norm, we
 let
$h=(v_{\theta},v_{\varphi})=v_{\theta}e_{\theta}+v_{\varphi}e_{\varphi}$.

The horizontal divergence $\hbox{div}$, the horizontal gradient
$\nabla=\hbox{grad}$, the horizontal covariant derivative
$\nabla_v$ and horizontal Laplace-Beltrami operator $\triangle$
for scalar and vector functions are defined by
$$\hbox{div}v=\hbox{div}(v_{\theta}e_{\theta}+v_{\varphi}e_{\varphi})=
\frac{1}{\sin \theta}(\frac{\partial
v_{\theta}\sin\theta}{\partial \theta}+\frac{\partial
v_{\varphi}}{\partial
\varphi}),\quad\quad\quad\quad\quad\quad\quad\quad\quad\eqno{(4.1)}$$
$$\nabla T=\hbox{grad}T=\frac{\partial
T}{\partial \theta}e_{\theta}+\frac{1}{\sin \theta}\frac{\partial
T}{\partial \varphi}e_{\varphi},
\quad\quad\quad\quad\quad\quad\quad\quad\quad\quad\quad\quad
\quad\quad\eqno{(4.2)}$$
$$\hbox{grad}\Phi_s=\frac{\partial
\Phi_s}{\partial \theta}e_{\theta}+\frac{1}{\sin
\theta}\frac{\partial \Phi_s}{\partial \varphi}e_{\varphi},
\quad\quad\quad\quad\quad\quad\quad\quad\quad\quad\quad\quad\quad
\quad\quad\eqno{(4.3)}$$
$$\nabla_v\widetilde{v}=(v_{\theta}\frac{\partial \widetilde{v_{\theta}}}{\partial
\theta}+\frac{v_{\varphi}}{\sin \theta}\frac{\partial
\widetilde{v_{\theta}}}{\partial
\varphi}-v_{\varphi}\widetilde{v_{\varphi}}\cot\theta)e_{\theta}+(v_{\theta}\frac{\partial
\widetilde{v_{\varphi}}}{\partial \theta}+\frac{v_{\varphi}}{\sin
\theta}\frac{\partial \widetilde{v_{\varphi}}}{\partial
\varphi}+v_{\varphi}\widetilde{v_{\theta}}\cot\theta)e_{\varphi},
\eqno{(4.4)}$$
$$\nabla_vT=v_{\theta}\frac{\partial T}{\partial
\theta}+\frac{v_{\varphi}}{\sin \theta}\frac{\partial T}{\partial
\varphi}, \quad\quad\quad\quad\quad\quad\quad\quad\quad
\quad\quad\quad\quad\quad\quad\quad\quad\quad\quad\eqno{(4.5)}$$
$$\nabla_vq=v_{\theta}\frac{\partial q}{\partial
\theta}+\frac{v_{\varphi}}{\sin \theta}\frac{\partial q}{\partial
\varphi}, \quad\quad\quad\quad\quad\quad\quad\quad\quad
\quad\quad\quad\quad\quad\quad\quad\quad\quad\quad\eqno{(4.6)}$$
$$\triangle T=\hbox{div}(\hbox{grad}T)=\frac{1}{\sin \theta}[\frac{\partial}{\partial\theta}(\sin\theta\frac{\partial
T}{\partial \theta})+\frac{1}{\sin \theta}\frac{\partial^2
T}{\partial
\varphi^2}],\quad\quad\quad\quad\quad\quad\quad\eqno{(4.7)}$$
$$\triangle q=\hbox{div}(\hbox{grad}q)=\frac{1}{\sin \theta}[\frac{\partial}{\partial\theta}(\sin\theta\frac{\partial
q}{\partial \theta})+\frac{1}{\sin \theta}\frac{\partial^2
q}{\partial
\varphi^2}],\quad\quad\quad\quad\quad\quad\quad\quad\eqno{(4.8)}$$
$$\triangle v=(\triangle v_{\theta}-\frac{2\cos\theta }{\sin^2 \theta}\frac{\partial
v_{\varphi}}{\partial \varphi}-\frac{v_{\theta}}{\sin^2\theta
})e_{\theta} +(\triangle v_{\varphi}+\frac{2\cos\theta }{\sin^2
\theta}\frac{\partial v_{\theta}}{\partial
\varphi}-\frac{v_{\varphi}}{\sin^2\theta })e_{\varphi},
\eqno{(4.9)}$$ where
$v=v_{\theta}e_{\theta}+v_{\varphi}e_{\varphi},\
\widetilde{v}=\widetilde{v_{\theta}}e_{\theta}+\widetilde{v_{\varphi}}e_{\varphi}\in
C^{\infty}(T\Omega|TS^2)$, $T,\ q\in C^{\infty}(\Omega)$,
$\Phi_s\in C^{\infty}(S^2)$.

Now we can define our working spaces for the problem (IBVP). Let
$$\widetilde{\mathcal{V}_1}:=\{v;\ v\in C^{\infty}(T\Omega|TS^2),\ \frac{\partial v}{\partial \xi}|_{\xi=0}=0,
\ \frac{\partial v}{\partial \xi}|_{\xi=1}=0,\ \int^1_0\hbox{div}
v\ d\xi=0\},\quad\quad\quad\quad$$
$$\widetilde{\mathcal{V}_2}:=\{T;\ T\in C^{\infty}(\Omega),\ \frac{\partial T}{\partial \xi}|_{\xi=0}=0,
\ \frac{\partial T}{\partial \xi}|_{\xi=1}=-\alpha_s
T\},\quad\quad\quad\quad\quad\quad\quad\quad\quad\quad$$
$$\widetilde{\mathcal{V}_3}:=\{q;\ q\in C^{\infty}(\Omega),\ \frac{\partial q}{\partial \xi}|_{\xi=0}=0,
\ \frac{\partial q}{\partial \xi}|_{\xi=1}=-\beta_s
q\},\quad\quad\quad\quad\quad\quad\quad\quad\quad\quad$$
$$V_1=\hbox{the closure of}\ \widetilde{\mathcal{V}_1} \hbox{with respect to the norm}\ \|\cdot\|_1,
\quad\quad\quad\quad\quad\quad\quad\quad\quad\quad\quad$$
$$V_2=\hbox{the closure of}\ \widetilde{\mathcal{V}_2} \hbox{with respect to the norm}\ \|\cdot\|_1,
\quad\quad\quad\quad\quad\quad\quad\quad\quad\quad\quad$$
$$V_3=\hbox{the closure of}\ \widetilde{\mathcal{V}_3} \hbox{with respect to the norm}\ \|\cdot\|_1,
\quad\quad\quad\quad\quad\quad\quad\quad\quad\quad\quad$$
$$H_1=\hbox{the closure of}\ \widetilde{\mathcal{V}_1} \hbox{with respect to the norm}\ |\cdot|_2,
\quad\quad\quad\quad\quad\quad\quad\quad\quad\quad\quad$$
$$H_2=L^2(\Omega),\quad\quad\quad\quad\quad\quad\quad\quad\quad\quad\quad\quad\quad\quad\quad\quad
\quad\quad\quad\quad\quad\quad\quad\quad\quad\quad\quad\quad\quad$$
$$V=V_1\times V_2\times V_3,\quad\quad\quad\quad\quad\quad\quad\quad\quad\quad\quad\quad
\quad\quad\quad\quad\quad\quad\quad\quad\quad\quad\quad\quad\quad\quad\quad$$
$$H=H_1\times H_2\times H_2.\quad\quad\quad\quad\quad\quad\quad\quad\quad\quad\quad\quad
\quad\quad\quad\quad\quad\quad\quad\quad\quad\quad\quad\quad\quad\quad$$
The inner products and norms on $V_1,\ V_2,\ V_3$ are given by
$$(v,v_1)_{V_1}=\int_\Omega(\nabla_{e_{\theta}}v\cdot\nabla_{e_{\theta}}v_1+
\nabla_{e_{\varphi}}v\cdot\nabla_{e_{\varphi}}v_1+\frac{\partial
v}{\partial \xi}\frac{\partial v_1}{\partial \xi}+v\cdot
v_1),\quad$$
$$\|v\|=(v,v)_{V_1}^{\frac{1}{2}},\ \ \forall v,\ v_1\in V_1,
\quad\quad\quad\quad\quad\quad\quad\quad\quad\quad\quad\quad\quad
\quad\quad$$
$$(T,T_1)_{V_2}=\int_\Omega(\hbox{grad}T\cdot\hbox{grad}T_1+\frac{\partial
T}{\partial \xi}\frac{\partial T_1}{\partial \xi}+TT_1),
\quad\quad\quad\quad\quad\quad\quad $$
$$\|T\|=(T,T)_{V_2}^{\frac{1}{2}},\ \ \forall T,\ T_1\in V_2,
\quad\quad\quad\quad\quad\quad\quad\quad\quad\quad\quad\quad\quad\quad$$
$$(q,q_1)_{V_3}=\int_\Omega(\hbox{grad}q\cdot\hbox{grad}q_1+\frac{\partial
q}{\partial \xi}\frac{\partial q_1}{\partial \xi}+qq_1),
\quad\quad\quad\quad\quad\quad\quad\quad $$
$$\|q\|=(q,q)_{V_3}^{\frac{1}{2}},\ \ \forall q,\ q_1\in V_2,
\quad\quad\quad\quad\quad\quad\quad\quad\quad\quad\quad\quad\quad
\quad\quad$$
$$(U,U_1)_{H}=(v,v_1)+(T,T_1)+(q,q_1),
\quad\quad\quad\quad\quad\quad\quad\quad\quad\quad\quad\quad $$
$$(U,U_1)_{V}=(v,v_1)_{V_1}+(T,T_1)_{V_2}+(q,q_1)_{V_3},\quad\quad\quad\quad
\quad\quad\quad\quad\quad$$
$$\quad\quad\quad\|U\|=(U,U)_{V}^{\frac{1}{2}},\ |U|_2=(U,U)_{H}^{\frac{1}{2}},
\ \forall U=(v,T,q),\ U_1=(v_1,T_1,q_1)\in V,$$ where
$(\cdot,\cdot)$ denotes the $L^2$ inner products in $H_1,\ H_2$.

\smallskip
\subsection{Some Lemmas}
\noindent {\bf Lemma 4.1}\quad Let $u=(u_{\theta}, u_{\varphi})$,
$u_1=((u_1)_{\theta}, (u_1)_{\varphi})\in
C^{\infty}(T\Omega|TS^2)$, and $p\in C^{\infty}(S^2)$. Then

(1) $$\displaystyle \int_{S^2}p  \hbox{div} u=-\displaystyle
\int_{S^2}\nabla p\cdot u , \eqno{(4.10)}$$ in particular,
$$\displaystyle \int_{S^2}\hbox{div} u=0;\eqno{(4.11)}$$

(2) $$\displaystyle \int_{\Omega}(-\triangle u)\cdot
u_1=\displaystyle
\int_{\Omega}(\nabla_{e_{\theta}}u\cdot\nabla_{e_{\theta}}u_1
+\nabla_{e_{\varphi}}u\cdot\nabla_{e_{\varphi}}u_1+u\cdot
u_1).\eqno{(4.12)}$$

\smallskip
\noindent{\bf Proof.} We can prove the first part of Lemma 4.1 by
using (4.1), (4.2) and Stokes Theorem(cf., e.g., \cite{T1, W1}).
From (4.4) and (4.9), by direct computation, we can obtain the
second part.

\smallskip
\noindent {\bf Lemma 4.2 (Interpolation Inequality)}\quad Let
$\Omega_1$ be a bounded domain in $\mathbb{R}^n$, whose boundary
$\partial \Omega_1$ satisfies $\partial \Omega_1\in C^m$. Then for
every $u\in W^{m, r}(\Omega_1)\cap L^q(\Omega_1)$, $0\leq l\leq
m$, one has

\smallskip
$1$)\ when $m-l-\frac{n}{r}$ is not a non-negative integer,
$$\|D^{l}u\|_{L^p(\Omega_1)}\leq c\|u\|_{W^{m,
r}(\Omega_1)}^{\alpha}\|u\|_{L^{q}(\Omega_1)}^{1-\alpha}, \ \
\hbox{where}\ l, p, \alpha, m, r, q \ \hbox{ satisfy}$$
$$\frac{l}{m}\leq \alpha\leq 1, 1\leq r, q\leq\infty, \frac{1}{p}-\frac{l}{n}
=\alpha(\frac{1}{r}-\frac{m}{n})+(1-\alpha)\frac{1}{q};$$

\smallskip
$2$)\ when $m-l-\frac{n}{r}$ is a non-negative integer,
$$\|D^{l}u\|_{L^p(\Omega_1)}\leq c\|u\|_{W^{m, r}(\Omega_1)}^{\alpha}\|u\|_{L^{q}(\Omega_1)}^{1-\alpha},
\ \ \hbox{where}\ l, p, \alpha, m, r, q \ \hbox{ satisfy}$$
$$\frac{l}{m}\leq \alpha<1, 1<r<\infty, 1<q<\infty.$$

In particular,

\smallskip
$i$)\ for $u\in H^1(S^2)$ (for the definitions of $H^1(S^2),\
L^p(S^2)$, cf. \cite{LM}),
$$\|u\|_{L^4(S^2)}\leq c\|u\|_{L^2(S^2)}^{\frac{1}{2}}\|u\|_{H^1(S^2)}^{\frac{1}{2}},\eqno{(4.13)}$$
$$\|u\|_{L^6(S^2)}\leq c\|u\|_{L^4(S^2)}^{\frac{2}{3}}\|u\|_{H^1(S^2)}^{\frac{1}{3}},\eqno{(4.14)}$$
$$\|u\|_{L^8(S^2)}\leq c\|u\|_{L^4(S^2)}^{\frac{1}{2}}\|u\|_{H^1(S^2)}^{\frac{1}{2}};\eqno{(4.15)}$$

\smallskip
$ii$)\ for $u\in H^1(\Omega)$,
$$\|u\|_{L^4(\Omega)}\leq c\|u\|_{L^2(\Omega)}^{\frac{1}{4}}\|u\|_{H^1(\Omega)}^{\frac{3}{4}}.\eqno{(4.16)}$$

\smallskip
\noindent{\bf Proof.} The proof of $1),\ 2)$ is similar to that of
the case $\Omega_1=\mathbb{R}^n$ but one must use the Extension
Theorem (for the detail of proof, we refer the reader to see,
e.g., \cite{A, G, T}). For the proof of $(4.13),\ (4.14),\
(4.15)$, one can refer to \cite[Chapter 1]{LM}.

\smallskip
\noindent {\bf Lemma 4.3}\quad For any $h\in C^{\infty}(S^2)$,
$v\in C^{\infty}(T\Omega|TS^2)$, we have
$$\displaystyle \int_{S^2}\nabla_{v}h+\displaystyle \int_{S^2}h \hbox{div}v=\displaystyle \int_{S^2}\hbox{div}(hv)=0.$$

\smallskip
\noindent{\bf Proof.} From (4.1), (4.5), (4.11), by direct
computation, we can prove Lemma 4.3.

\smallskip
\noindent {\bf Lemma 4.4}\quad Let $v$, $v_1\in V_1$, $T\in V_2,
q\in V_3$. Then we have

\smallskip
$1)$\ $\displaystyle \int_{\Omega}(\nabla_v
v_1+(\displaystyle\int_{\xi}^{1}\hbox{div}vd\xi^{'})\frac{\partial
v_1}{\partial \xi})v_1=0,$

\smallskip
$2)$\ $\displaystyle \int_{\Omega}(\nabla_v
T+(\displaystyle\int_{\xi}^{1}\hbox{div}vd\xi^{'})\frac{\partial
T}{\partial \xi})T=0,$

\smallskip
$3)$\ $\displaystyle \int_{\Omega}(\nabla_v
q+(\displaystyle\int_{\xi}^{1}\hbox{div}vd\xi^{'})\frac{\partial
q}{\partial \xi})q=0,$

\smallskip
$4)$\ $\displaystyle
\int_{\Omega}(\int_{\xi}^{1}\frac{bP}{p}\hbox{grad}[(1+aq)T]
d\xi^{'}\cdot v-\frac{bP}{p}(1+aq)W(v)\cdot T)=0.$

For the detail of proof for Lemma 4.4, we refer the reader to
\cite[Lemma 3.2]{GH}.

\smallskip
\noindent {\bf Lemma 4.5 (Minkowski Inequality)}\quad Let
$(X,\mu),\ (Y,\nu)$ be two measure spaces and $f(x,y)$ be a
measurable function about $\mu\times\nu$ on $X\times Y$. If for
a.e. $y\in Y,\ f(\cdot,y)\in L^p(X,\mu),\ 1\leq p\leq \infty,$ and
$\int_{Y}\|f(\cdot,y)\|_{L^p(X,\mu)}d\nu(y)< \infty$, then
$$\|\int_{Y}f(\cdot,y)d\nu(y)\|_{L^p(X,\mu)}\leq \int_{Y}\|f(\cdot,y)\|_{L^p(X,\mu)}d\nu(y).$$

\smallskip
\noindent {\bf Lemma 4.6(The Uniform Gronwall Lemma)}\quad Let
$\phi,\ \psi,\ \varphi$ be three positive locally integrable
functions on $[t_0,+\infty)$ such that $\varphi'$ is locally
integrable on $[t_0,+\infty)$, and which satisfy
$$\frac{d\varphi }{dt}\leq \phi \varphi +\psi \quad \hbox{for}\ t\geq t_0,
\quad\quad\quad\quad\quad\quad\quad\quad\quad\quad\quad\quad\quad\quad
\quad\quad\quad\quad$$
$$\int^{t+r}_t\phi (s)ds\leq a_1, \ \int^{t+r}_t\psi (s)ds\leq
a_2, \ \int^{t+r}_t\varphi (s)ds\leq a_3 \quad \hbox{for}\ t\geq
t_0,$$ where $r,\ a_1,\ a_2,\ a_3$ are positive constants. Then
$$\varphi (t+r)\leq (\frac{a_3}{r}+a_2)\exp(a_1),\quad \forall t\geq t_0.$$
For the detail of proof for Lemma 4.6, we refer the reader to
\cite[p91]{T1}.

\section{Global existence of strong solutions}
\subsection{Local existence of strong solutions}
In this subsection we recall the local, in time, existence of
strong solutions of the 3-D viscous primitive equations of the
large-scale moist atmosphere.

\noindent{\bf Definition 5.1} \quad Let $U_0=(v_0,T_0,q_0)\in V$,
and let $\mathcal{T}$ be a fixed positive time. $U=(v,T,q)$ is
called a strong solution of the system (2.11)-(2.17) on the time
interval $[0,\mathcal{T}]$ if it satisfies (2.11)-(2.13) in weak
sense such that
$$v\in C([0,\mathcal{T}];V_1)\cap L^2(0,\mathcal{T};(H^2(\Omega))^2),$$
$$T\in C([0,\mathcal{T}];V_2)\cap L^2(0,\mathcal{T};H^2(\Omega)),$$
$$q\in C([0,\mathcal{T}];V_3)\cap L^2(0,\mathcal{T};H^2(\Omega)),$$
$$\frac{\partial v}{\partial t}\in L^1(0,\mathcal{T};(L^2(\Omega))^2),\quad\quad\quad\quad\quad\ $$
$$\frac{\partial T}{\partial t}\in L^1(0,\mathcal{T};L^2(\Omega)),\quad\quad\quad\quad\quad\quad$$
$$\frac{\partial q}{\partial t}\in L^1(0,\mathcal{T};L^2(\Omega)).\quad\quad\quad\quad\quad\quad$$

\smallskip
\noindent{\bf Remark 5.2} \quad Since the 3-D viscous primitive
equations of large-scale moist atmosphere have not the time
derivative term of the vertical velocity $\omega=W(v)$, we can not
prove $\frac{\partial v}{\partial t}\in
L^2(0,\mathcal{T};(L^2(\Omega))^2),\ \frac{\partial T}{\partial
t},\ \frac{\partial q}{\partial t}\in
L^2(0,\mathcal{T};L^2(\Omega)).$

\smallskip
\noindent{\bf Proposition 5.3} \quad Let $Q_1,\ Q_2\in
H^1(\Omega),\ U_0=(v_0,T_0,q_0)\in V$. Then there exists
$\mathcal{T}_*>0,\ \mathcal{T}_*=\mathcal{T}_*(\|U_0\|)$, and
there exists a strong solution $U$ of the system (2.11)-(2.17) on
the interval $[0,\mathcal{T}_*]$.

\smallskip
\noindent{\bf Proof.} The proof of Proposition 5.3 is similar to
that for the primitive equations of large-scale ocean given in the
papers \cite{GM, TZ}. So we omit the detail of proof here.

In order to prove the global existence of strong solutions to the
system (2.11)-(2.17), we should make a priori estimates about
$H^1$-norm of the local solution $U(t)$ obtained in Proposition
5.3, that is, we should show that if $\mathcal{T}_*< \infty$ then
the $H^1$-norm of the strong solution $U(t)$ is bounded over the
interval $[0,\mathcal{T}_*]$.

\subsection{A priori estimates about local strong solutions}
\smallskip
\noindent{\bf $L^2$ estimates about $v$, $T$, $q$.} \quad Choosing
$v$ as a test function in equation (2.11), we obtain
\begin{eqnarray*}
&
&\frac{1}{2}\frac{d|v|_2^2}{dt}+\frac{1}{Re_1}\displaystyle\int_{\Omega}(|\nabla_{e_{\theta}}v|^2
+|\nabla_{e_{\varphi}}v|^2+|v|^2)+\frac{1}{Re_2}\displaystyle\int_{\Omega}|\frac{\partial
v}{\partial \xi}|^2\\
&=&-\displaystyle\int_{\Omega}(\nabla_{v}v+W(v)\frac{\partial
v}{\partial \xi}+\frac{f}{R_0}k\times v+\hbox{grad} \Phi_s)\cdot
v\\&-&\displaystyle\int_{\Omega}[\int_{\xi}^{1}\frac{b
P}{p}\hbox{grad}((1+aq)T)d\xi^{'}]\cdot v. \qquad \qquad \qquad
\qquad \qquad \qquad\quad\quad  (5.1)
\end{eqnarray*}
By Lemma 4.4, we get
$$\displaystyle\int_{\Omega}(\nabla_{v}v+W(v)\frac{\partial
v}{\partial \xi})\cdot v=0.\eqno{(5.2)}$$ $(\frac{f}{R_0}k\times
v)\cdot v=0$ implies
$$\displaystyle\int_{\Omega}(\frac{f}{R_0}k\times v)\cdot v=0.\eqno{(5.3)}$$
By using integration by parts and (2.14), we get
$$\displaystyle\int_{\Omega}\hbox{grad} \Phi_s\cdot v=-\displaystyle\int_{\Omega}\Phi_s\hbox{div}v
=-\displaystyle\int_{S^2}\Phi_s(\int_{0}^{1}\hbox{div}vd\xi)=0.\eqno{(5.4)}$$
So, from (5.1)-(5.4), we obtain
\begin{eqnarray*}
&
&\frac{1}{2}\frac{d|v|_2^2}{dt}+\frac{1}{Re_1}\displaystyle\int_{\Omega}(|\nabla_{e_{\theta}}v|^2
+|\nabla_{e_{\varphi}}v|^2+|v|^2)+\frac{1}{Re_2}\displaystyle\int_{\Omega}|\frac{\partial
v}{\partial \xi}|^2\\
&=&-\displaystyle\int_{\Omega}[\int_{\xi}^{1}\frac{b
P}{p}\hbox{grad}((1+aq)T)d\xi^{'}]\cdot v. \qquad \qquad \qquad
\qquad \qquad \qquad \quad   (5.5)
\end{eqnarray*}

Taking the inner product of equation (2.12) with $T$ in
$L^2(\Omega)$, we obtain
\begin{eqnarray*}
& &\frac{1}{2}\frac{d|T|_2^2}{dt}+\frac{1}{R
t_1}\displaystyle\int_{\Omega}|\nabla T|^2+\frac{1}{R
t_2}\displaystyle\int_{\Omega}|\frac{\partial T}{\partial
\xi}|^2+\frac{\alpha_s}{R t_2}|T|_{\xi=1}|_2^2\\
&=&-\displaystyle\int_{\Omega}(\nabla_v T+W(v)\frac{\partial
T}{\partial \xi})T+\displaystyle\int_{\Omega}\frac{b P}{p}(1+a q)T
W(v)+\displaystyle\int_{\Omega}Q_1 T. \quad \quad \quad (5.6)
\end{eqnarray*}
By Lemma 4.4, we have
$$-\displaystyle\int_{\Omega}(\nabla_v T+W(v)\frac{\partial
T}{\partial \xi})T=0.   \eqno(5.7)$$ Combining (5.6) with (5.7),
we get

\begin{eqnarray*}
& &\frac{1}{2}\frac{d|T|_2^2}{dt}+\frac{1}{R
t_1}\displaystyle\int_{\Omega}|\nabla T|^2+\frac{1}{R
t_2}\displaystyle\int_{\Omega}|\frac{\partial T}{\partial
\xi}|^2+\frac{\alpha_s}{R t_2}|T|_{\xi=1}|_2^2\\
&=&\displaystyle\int_{\Omega}\frac{b P}{p}(1+a q)T
W(v)+\displaystyle\int_{\Omega}Q_1 T. \quad \quad \quad \qquad
\qquad \qquad \qquad \qquad \qquad (5.8)
\end{eqnarray*}
Similarly to (5.8), we have
$$\frac{1}{2}\frac{d|q|_2^2}{dt}+\frac{1}{R
q_1}\displaystyle\int_{\Omega}|\nabla q|^2+\frac{1}{R
q_2}\displaystyle\int_{\Omega}|\frac{\partial q}{\partial
\xi}|^2+\frac{\beta_s}{R
q_2}|q|_{\xi=1}|_2^2=\displaystyle\int_{\Omega}qQ_2 . \eqno(5.9)$$
By using Lemma 4.4, we have
$$-\displaystyle\int_{\Omega}[\int_{\xi}^{1}\frac{b
P}{p}\hbox{grad}((1+aq)T)d\xi^{'}]\cdot
v+\displaystyle\int_{\Omega}\frac{b P}{p}(1+cq)TW(v)=0.
\eqno(5.10)$$ From (5.5) and (5.8)-(5.10), we obtain
\begin{eqnarray*}
&
&\frac{1}{2}\frac{d(|v|_2^2+|T|_2^2+|q|_2^2)}{dt}+\frac{1}{Re_1}\displaystyle\int_{\Omega}(|\nabla_{e_{\theta}}v|^2
+|\nabla_{e_{\varphi}}v|^2+|v|^2)+\frac{1}{Re_2}\displaystyle\int_{\Omega}|\frac{\partial
v}{\partial \xi}|^2\\&+&\frac{1}{R
t_1}\displaystyle\int_{\Omega}|\nabla T|^2+\frac{1}{R
t_2}\displaystyle\int_{\Omega}|\frac{\partial T}{\partial
\xi}|^2+\frac{\alpha_s}{R t_2}|T|_{\xi=1}|_2^2+\frac{1}{R
q_1}\displaystyle\int_{\Omega}|\nabla q|^2\\&+&\frac{1}{R
q_2}\displaystyle\int_{\Omega}|\frac{\partial q}{\partial
\xi}|^2+\frac{\beta_s}{R q_2}|q|_{\xi=1}|_2^2\\
&=&\displaystyle\int_{\Omega}Q_1 T+\displaystyle\int_{\Omega}qQ_2.
\qquad \qquad \qquad \qquad \qquad \qquad \qquad \qquad \qquad
\qquad \  (5.11)
\end{eqnarray*}
By $T(\theta, \varphi,
\xi)=-\displaystyle\int_{\xi}^{1}\frac{\partial T}{\partial
\xi^{'}}d\xi^{'}+T|_{\xi=1}$, using H\"older inequality and
Cauchy-Schwarz inequality, we have
$$|T|_2^2\leq2|\frac{\partial T}{\partial
\xi}|_2^2+2|T|_{\xi=1}|_2^2. \eqno(5.12)$$  Similarly to (5.12),
we get
$$|q|_2^2\leq2|\frac{\partial q}{\partial
\xi}|_2^2+2|q|_{\xi=1}|_2^2. \eqno(5.13)$$ By Young inequality, we
have
$$|\displaystyle\int_{\Omega}Q_1 T|\leq
\varepsilon|T|_2^2+c|Q_1|_2^2,\quad
|\displaystyle\int_{\Omega}qQ_2 |\leq
\varepsilon|q|_2^2+c|Q_2|_2^2. \eqno(5.14)$$ In this article, $c$
will denote positive constant and can be determined in concrete
conditions. $\varepsilon$ is a small enough positive constant.
Therefore, we obtain
\begin{eqnarray*}
&
&\frac{d(|v|_2^2+|T|_2^2+|q|_2^2)}{dt}+\frac{1}{Re_1}\displaystyle\int_{\Omega}(|\nabla_{e_{\theta}}v|^2
+|\nabla_{e_{\varphi}}v|^2+|v|^2)+\frac{1}{Re_2}\displaystyle\int_{\Omega}|\frac{\partial
v}{\partial \xi}|^2\\&+&\frac{1}{R
t_1}\displaystyle\int_{\Omega}|\nabla T|^2+\frac{1}{R
t_2}\displaystyle\int_{\Omega}|\frac{\partial T}{\partial
\xi}|^2+\frac{\alpha_s}{R t_2}|T|_{\xi=1}|_2^2+\frac{1}{R
q_1}\displaystyle\int_{\Omega}|\nabla q|^2\\&+&\frac{1}{R
q_2}\displaystyle\int_{\Omega}|\frac{\partial q}{\partial
\xi}|^2+\frac{\beta_s}{R q_2}|q|_{\xi=1}|_2^2\\
&\leq&c\displaystyle\int_{\Omega}(Q_1^2+Q_2^2). \qquad \qquad
\qquad \qquad \qquad \qquad \qquad \qquad \qquad \quad \quad\quad
\ \ (5.15)
\end{eqnarray*}
By (5.12), (5.13), (5.15) and thanks to the Gronwall  inequality,
we have
$$|v|_2^2+|T|_2^2+|q|_2^2\leq e^{-c_0 t}(|v_0|_2^2+|T_0|_2^2+|q_0|_2^2)+c(|Q_1|_2^2+|Q_2|_2^2)
\leq E_0, \eqno(5.16)$$ where $c_0=\hbox{min}\{\frac{1}{Re_1},
\frac{1}{2Rt_2}, \frac{\alpha_s}{2R_{t_2}}, \frac{1}{2Rq_2},
\frac{\beta_s}{2Rq_2}\}>0,\ t\geq 0$ and $E_0$ is a positive
constant. By Minkowski inequality and H\"older inequality, for any
$t\geq 0$ we have
\begin{eqnarray*}
\|\bar{v}(t)\|_{L^2(S^2)}^2
&\leq&|v(t)|_2^2\\
&\leq&e^{-c_0
t}(|v_0|_2^2+|T_0|_2^2+|q_0|_2^2)+c(|Q_1|_2^2+|Q_2|_2^2)\leq E_0
.\quad\  (5.17)
\end{eqnarray*}
From (5.12), (5.13), (5.15) and (5.16), we get
\begin{eqnarray*}
&
&c_1\displaystyle\int_{t}^{t+r}[\int_{\Omega}(|\nabla_{e_{\theta}}v|^2
+|\nabla_{e_{\varphi}}v|^2+|\frac{\partial v}{\partial
\xi}|^2)+\int_{\Omega}(|\nabla T|^2+|\frac{\partial T}{\partial
\xi}|^2+|T|^2)\\
&+&\int_{\Omega}(|\nabla
q|^2+|\frac{\partial q}{\partial \xi}|^2+|q|^2)+|T|_{\xi=1}|_2^2+|q|_{\xi=1}|_2^2]+|U(t)|_2^2\\
&\leq&2e^{-c_0
t}(|v_0|_2^2+|T_0|_2^2+|q_0|_2^2)+c(|Q_1|_2^2+|Q_2|_2^2)(2+r)\leq
E_1,\ \quad\quad (5.18)
\end{eqnarray*}
where $c_1=\hbox{min}\{\frac{1}{Re_1}, \frac{1}{Re_2}, \frac{1}{R
t_1}, \frac{1}{3R t_2}, \frac{\alpha_s}{2R t_2}, \frac{1}{Rq_1},
\frac{1}{3Rq_2}, \frac{\beta_s}{2Rq_2}\}$, $t\geq 0$, $1\geq r>0$
given, $E_1$ is a positive constant and
$\displaystyle\int_{t}^{t+r}\cdot ds$ is denoted by
$\displaystyle\int_{t}^{t+r}\cdot$. Since
$$\int_{S^2}(|\nabla_{e_{\theta}}\bar{v}|^2
+|\nabla_{e_{\varphi}}\bar{v}|^2)\leq
\int_{\Omega}(|\nabla_{e_{\theta}}v|^2
+|\nabla_{e_{\varphi}}v|^2),$$ from (5.18), we have
$$c_1\displaystyle\int_{t}^{t+r}\int_{S^2}(|\nabla_{e_{\theta}}\bar{v}|^2
+|\nabla_{e_{\varphi}}\bar{v}|^2)+\| \bar{v}\|_{L^2(S^2)}^2\leq
E_1, \forall t\geq 0.\eqno{(5.19)}$$

\noindent{\bf $L^4$ estimates about $q$} \quad By taking the inner
product of equation (2.13) with $|q|^2q$ in $L^2(\Omega)$, we get
\begin{eqnarray*}
& &\frac{1}{4}\frac{d
|q|_4^4}{dt}+\frac{3}{Rq_1}\displaystyle\int_{\Omega}|\nabla
q|^2q^2+\frac{3}{Rq_2}\displaystyle\int_{\Omega}|\frac{\partial
q}{\partial
\xi}|^2q^2+\frac{\beta_s}{Rq_2}\displaystyle\int_{S^2}|q|_{\xi=1}|^4\\
&=&\displaystyle\int_{\Omega}Q_2|q|^2q-\displaystyle\int_{\Omega}(\nabla_{v}q+(\int_{\xi}^{1}\hbox{div}
vd\xi^{'})\frac{\partial q}{\partial \xi})|q|^2q. \qquad \qquad
\qquad \qquad \ (5.20)
\end{eqnarray*}
By Lemma 4.3, we have
\begin{eqnarray*}
&
&\displaystyle\int_{\Omega}(\nabla_{v}q+(\int_{\xi}^{1}\hbox{div}
vd\xi^{'})\frac{\partial q}{\partial \xi})|q|^2q\\&=&
\frac{1}{4}\displaystyle\int_{\Omega}\nabla_{v}q^4+\displaystyle\int_{S^2}[\int_{0}^{1}(\int_{\xi}^{1}\hbox{div}
vd\xi^{'})d(\frac{1}{4}q^4)]\\
&=&\frac{1}{4}\displaystyle\int_{\Omega}(\nabla_{v}q^4+q^4\hbox{div}v)\\
&=&0.\qquad
\qquad\qquad\qquad\qquad\qquad\qquad\qquad\qquad\qquad\qquad\qquad\qquad\quad(5.21)
\end{eqnarray*}
Combining (5.20) with (5.21), we obtain $$\frac{1}{4}\frac{d
|q|_4^4}{dt}+\frac{3}{Rq_1}\displaystyle\int_{\Omega}|\nabla
q|^2q^2+\frac{3}{Rq_2}\displaystyle\int_{\Omega}|\frac{\partial
q}{\partial
\xi}|^2q^2+\frac{\beta_s}{Rq_2}\displaystyle\int_{S^2}|q|_{\xi=1}|^4
=\displaystyle\int_{\Omega}Q_2|q|^2q.\eqno{(5.22)}
$$
Since $q^4(\theta, \varphi,
\xi)=-\displaystyle\int_{\xi}^{1}\frac{\partial q^4}{\partial
\xi^{'}}d\xi^{'}+q^4|_{\xi=1}$, we get by using H\"older
inequality and Cauchy-Schwarz inequality,
\begin{eqnarray*}
|q|_4^4&\leq&4\displaystyle\int_{S^2}(\int_{0}^{1}(\int_{\xi}^{1}|q|^3|\frac{\partial
q}{\partial \xi^{'}}|d\xi^{'})d\xi)+|q|_{\xi=1}|_4^4\\
&\leq&c\displaystyle\int_{S^2}[(\int_{0}^{1}|q|^2|\frac{\partial
q}{\partial
\xi}|^2d\xi)^{\frac{1}{2}}(\int_{0}^{1}q^4d\xi)^{\frac{1}{2}}]+|q|_{\xi=1}|_4^4\\
&\leq&c(\displaystyle\int_{\Omega}|q|^2|\frac{\partial q}{\partial
\xi}|^2)+\frac{1}{2}\displaystyle\int_{\Omega}q^4+|q|_{\xi=1}|_4^4.\qquad\qquad\qquad\qquad\qquad\
\ \ (5.23)
\end{eqnarray*}
Since $$|\displaystyle\int_{\Omega}Q_2|q|^2q|\leq
c|Q_2|_4^4+\varepsilon|q|_4^4,$$ from (5.22)-(5.23), by choosing
$\varepsilon$ small enough,  we obtain
$$\frac{d|q|_4^4}{d t}+c_2|q|_4^4\leq c|Q_2|_4^4,\eqno(5.24)$$
where $c_2$ is a positive constant. By Gronwall inequality, we
have
$$|q(t)|_4^4\leq e^{-c_2 t}|q_0|_4^4+c|Q_2|_4^4\leq E_2,\eqno(5.25)$$
where $t\geq 0,\ E_2$ is a positive constant. From (5.22) and
(5.25), we get
$$c_1\int_t^{t+r}|q|_{\xi=1}|^4_4\leq 2E_2, \ \hbox{for any}\ t\geq 0.\eqno(5.26)$$

\noindent{\bf $L^3$ estimates about $T$} \quad We take the inner
product of equation (2.12) with $|T|T$ in $L^2(\Omega)$ and obtain
\begin{eqnarray*}
& &\frac{1}{3}\frac{d
|T|_3^3}{dt}+\frac{2}{Rt_1}\displaystyle\int_{\Omega}|\nabla
T|^2|T|+\frac{2}{Rt_2}\displaystyle\int_{\Omega}|\frac{\partial
T}{\partial
\xi}|^2|T|+\frac{\alpha_s}{Rt_2}\displaystyle\int_{S^2}|T|_{\xi=1}|^3\\
&=&\displaystyle\int_{\Omega}Q_1|T|T-\displaystyle\int_{\Omega}(\nabla_v
T+(\int_{\xi}^{1}\hbox{div} vd\xi^{'})\frac{\partial T}{\partial
\xi})|T|T\\
&+&\displaystyle\int_{\Omega}\frac{b
P}{p}(\int_{\xi}^{1}\hbox{div}
vd\xi^{'})|T|T+\displaystyle\int_{\Omega}\frac{ab
P}{p}(\int_{\xi}^{1}\hbox{div}vd\xi^{'})q|T|T.\qquad \qquad \
(5.27)
\end{eqnarray*}
Using H\"older inequality and Young inequality, we have
$$|\displaystyle\int_{\Omega}Q_1|T|T|\leq c|Q_1|_3^3+\varepsilon|T|_3^3.\eqno(5.28)$$
By H\"older inequality, Minkowski inequality and (4.14), we get
\begin{eqnarray*}
&&|\displaystyle\int_{\Omega}\frac{b
P}{p}(\int_{\xi}^{1}\hbox{div}
vd\xi^{'})|T|T|\\
&\leq&c\int_{0}^{1}[(\displaystyle\int_{S^2}(\int_{\xi}^{1}\hbox{div}
vd\xi^{'})^2)^{\frac{1}{2}}(\displaystyle\int_{S^2}|T|^4)^{\frac{1}{2}}]d\xi\\
&\leq&c\|(\displaystyle\int_{S^2}(\int_{\xi}^{1}\hbox{div}
vd\xi^{'})^2)^{\frac{1}{2}}\|_{L_{\xi}^{\infty}}\int_{0}^{1}\|T\|_{L^2(S^2)}\|T\|_{H^1(S^2)}d\xi\\
&\leq&c(\displaystyle\int_{\Omega}(|\nabla_{e_{\varphi}}v|^2
+|\nabla_{e_{\varphi}}v|^2))^{\frac{1}{2}}|T|_2\|T\|\\
&\leq&c\displaystyle\int_{\Omega}(|\nabla_{e_{\varphi}}v|^2
+|\nabla_{e_{\varphi}}v|^2)+c|T|_2^2\|T\|^2.\
\quad\quad\quad\quad\quad\quad\quad\quad\quad\quad\quad\quad(5.29)
\end{eqnarray*}
By H\"older inequality, Young inequality and
$$\|u\|_{L^{\frac{16}{3}}(S^2)}\leq c\|u\|^{\frac{3}{8}}_{L^2(S^2)}
\|u\|^{\frac{5}{8}}_{H^1(S^2)},\ \hbox{for any}\ u\in H^1(S^2),$$
we have
\begin{eqnarray*}
& &|\displaystyle\int_{\Omega}\frac{ab
P}{p}q(\int_{\xi}^{1}\hbox{div}
vd\xi^{'})|T|T|\\
&\leq&c|\int_{0}^{1}[\displaystyle\int_{S^2}q(\int_{\xi}^{1}\hbox{div}
vd\xi^{'})|T|T]d\xi|\\
&\leq&c\displaystyle\int_{0}^{1}[(\displaystyle\int_{S^2}q^4)^{\frac{1}{4}}
(\displaystyle\int_{S^2}(|T|^\frac{3}{2})^\frac{16}{3})^{\frac{1}{4}}(\displaystyle\int_{S^2}(\int_{\xi}^{1}\hbox{div}
vd\xi^{'})^2)^{\frac{1}{2}}]d\xi\\
&\leq&c\displaystyle\int_{0}^{1}[(\displaystyle\int_{S^2}q^4)^{\frac{1}{4}}
(\||T|^\frac{3}{2}\|_{L^2(S^2)}^2\||T|^\frac{3}{2}\|_{H^1(S^2)}^\frac{10}{3})
^{\frac{1}{4}}(\displaystyle\int_{S^2}(\int_{\xi}^{1}\hbox{div}
vd\xi^{'})^2)^{\frac{1}{2}}]d\xi\\
&\leq&c\displaystyle\int_{0}^{1}[(\displaystyle\int_{S^2}q^4)^{\frac{1}{4}}\|T\|_{L^3(S^2)}^{\frac{3}{4}}
\||T|^\frac{3}{2}\|_{H^1(S^2)}^{\frac{5}{6}}(\displaystyle\int_{S^2}(\int_{\xi}^{1}\hbox{div}
vd\xi^{'})^2)^{\frac{1}{2}}]d\xi\\
&\leq&c|q|_4|T|_3^{\frac{3}{4}}(\displaystyle\int_{\Omega}|T||\nabla
T|^2+|T|_3^3)^{\frac{5}{12}}\|(\displaystyle\int_{S^2}(\int_{\xi}^{1}\hbox{div}
vd\xi^{'})^2)^{\frac{1}{2}}\|_{L_{\xi}^{\infty}}\\
&\leq&c|q|_4|T|_3^{\frac{3}{4}}(\displaystyle\int_{\Omega}|T||\nabla
T|^2+|T|_3^3)^{\frac{5}{12}}\|\displaystyle\int_{S^2}\int_{\xi}^{1}
(|\nabla_{e_{\theta}}v|^2+|\nabla_{e_{\varphi}}v|^2)d\xi^{'}\|_{L_{\xi}^{\infty}}^{\frac{1}{2}}\\
&\leq&c|q|_4|T|_3^{\frac{3}{4}}(\displaystyle\int_{\Omega}|T||\nabla
T|^2+|T|_3^3)^{\frac{5}{12}}[\displaystyle\int_{\Omega}(|\nabla_{e_{\theta}}v|^2
+|\nabla_{e_{\varphi}}v|^2)]^{\frac{1}{2}}\\
&\leq&c|q|_4^{\frac{12}{7}}|T|_3^{\frac{9}{7}}[\displaystyle\int_{\Omega}(|\nabla_{e_{\theta}}v|^2
+|\nabla_{e_{\varphi}}v|^2)]^{\frac{6}{7}}+\varepsilon(\displaystyle\int_{\Omega}|T||\nabla
T|^2+|T|_3^3)\\
&\leq&c|q|_4^{\frac{12}{7}}(1+|T|_3^3)[1+\displaystyle\int_{\Omega}(|\nabla_{e_{\theta}}v|^2
+|\nabla_{e_{\varphi}}v|^2)]+\varepsilon(\displaystyle\int_{\Omega}|T||\nabla
T|^2+|T|_3^3). \   (5.30)
\end{eqnarray*}
Choosing $\varepsilon$ small enough and using a inequality similar
to (5.23), we derive from (5.27)-(5.30)
\begin{eqnarray*}
& &\frac{d
|T|_3^3}{dt}+\frac{2}{Rt_1}\displaystyle\int_{\Omega}|\nabla
T|^2|T|+\frac{2}{Rt_2}\displaystyle\int_{\Omega}|\frac{\partial
T}{\partial
\xi}|^2|T|+\frac{\alpha_s}{Rt_2}\displaystyle\int_{S^2}|T|_{\xi=1}|^3\\
&\leq&c|q|_4^{\frac{12}{7}}[1+\displaystyle\int_{\Omega}(|\nabla_{e_{\theta}}v|^2
+|\nabla_{e_{\varphi}}v|^2)]|T|_3^3+c|Q_1|_3^3\\
&+&c(1+|q|_4^{\frac{12}{7}})[1+\displaystyle\int_{\Omega}(|\nabla_{e_{\theta}}v|^2
+|\nabla_{e_{\varphi}}v|^2)]+c|T|_2^2\|T\|^2.
\quad\quad\quad\quad\quad\quad (5.31)
\end{eqnarray*}
By Lemma 4.6, (5.18), (5.25) and $|T|_3^3\leq
c|T|_2^{\frac{3}{2}}\|T\|^{\frac{3}{2}}$, we obtain
$$|T(t+r)|_3^3\leq c(|Q_1|_3^3+(1+E_2^{\frac{3}{7}})(1+E_1)+E_1^2+\frac{(E_0E_1)^{\frac{3}{4}}}{r})
\exp(cE_2^{\frac{3}{7}}(1+E_1))= E_3,\eqno{(5.32)}$$ where $E_3$
is a positive constant, $t\geq 0$. By Gronwall inequality, we
prove
$$|T(t)|_3^3\leq c(|Q_1|_3^3+(1+E_2^{\frac{3}{7}})(1+E_1)+|T_0|_3^3)
\exp(cE_2^{\frac{3}{7}}(1+E_1)),\ \hbox{for any}\ 0\leq
t<r.\eqno{(5.33)}$$

\noindent{\bf $L^4$ estimates about $T$} \quad We take the inner
product of equation (2.12) with $|T|^2T$ in $L^2(\Omega)$ and
obtain
\begin{eqnarray*}
& &\frac{1}{4}\frac{d
|T|_4^4}{dt}+\frac{3}{Rt_1}\displaystyle\int_{\Omega}|\nabla
T|^2T^2+\frac{3}{Rt_2}\displaystyle\int_{\Omega}|\frac{\partial
T}{\partial
\xi}|^2T^2+\frac{\alpha_s}{Rt_2}\displaystyle\int_{S^2}|T|_{\xi=1}|^4\\
&=&\displaystyle\int_{\Omega}Q_1|T|^2T-\displaystyle\int_{\Omega}(\nabla_v
T+(\int_{\xi}^{1}\hbox{div} vd\xi^{'})\frac{\partial T}{\partial
\xi})|T|^2T\\
&+&\displaystyle\int_{\Omega}\frac{b
P}{p}(\int_{\xi}^{1}\hbox{div}
vd\xi^{'})|T|^2T+\displaystyle\int_{\Omega}\frac{ab
P}{p}(\int_{\xi}^{1}\hbox{div}vd\xi^{'})q|T|^2T.\qquad \qquad \
(5.34)
\end{eqnarray*}
Using H\"older inequality and Young inequality, we have
$$|\displaystyle\int_{\Omega}Q_1|T|^2T|\leq c|Q_1|_4^4+\varepsilon|T|_4^4.\eqno(5.35)$$
By H\"older inequality, Minkowski inequality and (4.14), we get
\begin{eqnarray*}
&&|\displaystyle\int_{\Omega}\frac{b
P}{p}(\int_{\xi}^{1}\hbox{div}
vd\xi^{'})|T|^2T|\\
&\leq&c\int_{0}^{1}[(\displaystyle\int_{S^2}(\int_{\xi}^{1}\hbox{div}
vd\xi^{'})^2)^{\frac{1}{2}}(\displaystyle\int_{S^2}|T|^6)^{\frac{1}{2}}]d\xi\\
&\leq&c\|(\displaystyle\int_{S^2}(\int_{\xi}^{1}\hbox{div}
vd\xi^{'})^2)^{\frac{1}{2}}\|_{L_{\xi}^{\infty}}\int_{0}^{1}\|T\|_{L^4(S^2)}^2\|T\|_{H^1(S^2)}d\xi\\
&\leq&c(\displaystyle\int_{\Omega}(|\nabla_{e_{\varphi}}v|^2
+|\nabla_{e_{\varphi}}v|^2))^{\frac{1}{2}}|T|_4^2\|T\|\\
&\leq&c\displaystyle\int_{\Omega}(|\nabla_{e_{\varphi}}v|^2
+|\nabla_{e_{\varphi}}v|^2)+c|T|_4^4\|T\|^2.\
\quad\quad\quad\quad\quad\quad\quad\quad\quad\quad\quad\quad\
(5.36)
\end{eqnarray*}
By H\"older inequality, Young inequality and
$$\|u\|_{L^6(S^2)}\leq c\|u\|^{\frac{1}{3}}_{L^2(S^2)}
\|u\|^{\frac{2}{3}}_{H^1(S^2)},\ \hbox{for any}\ u\in H^1(S^2),$$
we have
\begin{eqnarray*}
& &|\displaystyle\int_{\Omega}\frac{ab
P}{p}q(\int_{\xi}^{1}\hbox{div}
vd\xi^{'})|T|^2T|\\
&\leq&c|\int_{0}^{1}[\displaystyle\int_{S^2}|q(\int_{\xi}^{1}\hbox{div}
vd\xi^{'})|T|^2T|]d\xi|\\
&\leq&c\displaystyle\int_{0}^{1}[(\displaystyle\int_{S^2}q^4)^{\frac{1}{4}}
(\displaystyle\int_{S^2}(|T|^2)^6)^{\frac{1}{4}}(\displaystyle\int_{S^2}(\int_{\xi}^{1}\hbox{div}
vd\xi^{'})^2)^{\frac{1}{2}}]d\xi\\
&\leq&c\displaystyle\int_{0}^{1}[(\displaystyle\int_{S^2}q^4)^{\frac{1}{4}}
(\|T^2\|_{L^2(S^2)}^2\|T^2\|_{H^1(S^2)}^4)^{\frac{1}{4}}(\displaystyle\int_{S^2}(\int_{\xi}^{1}\hbox{div}
vd\xi^{'})^2)^{\frac{1}{2}}]d\xi\\
&\leq&c\displaystyle\int_{0}^{1}[(\displaystyle\int_{S^2}q^4)^{\frac{1}{4}}\|T\|_{L^4(S^2)}
\|T^2\|_{H^1(S^2)}(\displaystyle\int_{S^2}(\int_{\xi}^{1}\hbox{div}
vd\xi^{'})^2)^{\frac{1}{2}}]d\xi\\
&\leq&c|q|_4|T|_4(\displaystyle\int_{\Omega}|T|^2|\nabla
T|^2+|T|_4^4)^{\frac{1}{2}}\|(\displaystyle\int_{S^2}(\int_{\xi}^{1}\hbox{div}
vd\xi^{'})^2)^{\frac{1}{2}}\|_{L_{\xi}^{\infty}}\\
&\leq&c|q|_4|T|_4(\displaystyle\int_{\Omega}|T|^2|\nabla
T|^2+|T|_4^4)^{\frac{1}{2}}[\displaystyle\int_{\Omega}(|\nabla_{e_{\theta}}v|^2
+|\nabla_{e_{\varphi}}v|^2)]^{\frac{1}{2}}\\
&\leq&c|q|_4^2|T|_4^2\displaystyle\int_{\Omega}(|\nabla_{e_{\theta}}v|^2
+|\nabla_{e_{\varphi}}v|^2)+\varepsilon(\displaystyle\int_{\Omega}|T|^2|\nabla
T|^2+|T|_4^4)\\
&\leq&c(|q|_4^4+|T|_4^4)\displaystyle\int_{\Omega}(|\nabla_{e_{\theta}}v|^2
+|\nabla_{e_{\varphi}}v|^2)+\varepsilon(\displaystyle\int_{\Omega}|T|^2|\nabla
T|^2+|T|_4^4). \quad \quad  (5.37)
\end{eqnarray*}
Choosing $\varepsilon$ small enough and using a inequality similar
to (5.23), we derive from (5.34)-(5.37),
\begin{eqnarray*}
& &\frac{d
|T|_4^4}{dt}+\frac{3}{Rt_1}\displaystyle\int_{\Omega}|\nabla
T|^2T^2+\frac{3}{Rt_2}\displaystyle\int_{\Omega}|\frac{\partial
T}{\partial
\xi}|^2T^2+\frac{\alpha_s}{Rt_2}\displaystyle\int_{S^2}|T|_{\xi=1}|^4\\
&\leq&c[\|T\|^2+\displaystyle\int_{\Omega}(|\nabla_{e_{\theta}}v|^2
+|\nabla_{e_{\varphi}}v|^2)]|T|_4^4+c|Q_1|_4^4\\
&+&c(1+|q|_4^4)\displaystyle\int_{\Omega}(|\nabla_{e_{\theta}}v|^2
+|\nabla_{e_{\varphi}}v|^2).
\quad\quad\quad\quad\quad\quad\quad\quad\quad\qquad  \qquad \ \
(5.38)
\end{eqnarray*}
By Lemma 4.6, (5.18), (5.25), (5.32) and $|T|_4^4\leq
c|T|_3^2\|T\|^2$, we obtain
$$|T(t+2r)|_4^4\leq c(|Q_1|_4^4+E_1E_2+E_1+\frac{E_1E_3^{\frac{2}{3}}}{r})\exp(cE_1)=
E_4,\eqno{(5.39)}$$ where $E_4$ is a positive constant and $t\geq
0$. By Gronwall inequality, from (5.38) we prove
$$|T(t)|_4^4\leq c(|Q_1|_4^4+E_1E_2+E_1+|T_0|_4^4)\exp(cE_1)=C_1,
$$ where $C_1=C_1(\|U_0\|,\ \|Q_1\|_1,\  \|Q_2\|_1)>0$ and $0\leq t<2r$. From (5.38) and (5.39), we get
$$c_1\int_{t+2r}^{t+3r}|T|_{\xi=1}|^4_4\leq E_4^2+E_4, \ \hbox{for any}\ t\geq 0.\eqno(5.40)$$

\noindent{\bf $L^3$ estimates about $\tilde{v}$} \quad We take the
inner product of equation (2.23) with $|\tilde{v}|\tilde{v}$ in
$L^2(\Omega)\times L^2(\Omega)$, and obtain
\begin{eqnarray*}
& &\frac{1}{3}\frac{d
|\tilde{v}|_3^3}{dt}+\frac{1}{Re_1}\displaystyle\int_{\Omega}[(|\nabla_{e_{\theta}}\tilde{v}|^2
+|\nabla_{e_{\varphi}}\tilde{v}|^2)|\tilde{v}|+\frac{4}{9}|\nabla_{e_{\theta}}|\tilde{v}|^{\frac{3}{2}}|^2
+\frac{4}{9}|\nabla_{e_{\varphi}}|\tilde{v}|^{\frac{3}{2}}|^2+|\tilde{v}|^3]\\
&+&\frac{1}{Re_2}\displaystyle\int_{\Omega}(|\tilde{v}_\xi|^2|\tilde{v}|
+\frac{4}{9}|\partial_\xi|\tilde{v}|^{\frac{3}{2}}|^2)\\
&=&-\displaystyle\int_{\Omega}(\nabla_{\tilde{v}}\tilde{v}+(\int_{\xi}^{1}\hbox{div}\tilde{v}d\xi^{'})\frac{\partial
\tilde{v}}{\partial \xi})\cdot
|\tilde{v}|\tilde{v}-\displaystyle\int_{\Omega}(\nabla_{\bar{v}}\tilde{v})\cdot|\tilde{v}|\tilde{v}
-\displaystyle\int_{\Omega}|\tilde{v}|\tilde{v}\cdot\nabla_{\tilde{v}}\bar{v}\\
&-&\displaystyle\int_{\Omega}[\int_{\xi}^{1}\frac{bP}{p}\hbox{grad}((1+aq)T)d\xi^{'}
-\int_{0}^{1}(\int_{\xi}^{1}\frac{bP}{p}\hbox{grad}((1+aq)T)d\xi^{'})d\xi]\cdot|\tilde{v}|\tilde{v}\\
&+& \displaystyle\int_{\Omega}\overline{(\tilde{v}\hbox{div}
\tilde{v}+\nabla_{\tilde{v}}\tilde{v})}\cdot|\tilde{v}|\tilde{v}-\displaystyle\int_{\Omega}(\frac{f}{R_0}k
\times\tilde{v})\cdot|\tilde{v}|\tilde{v}, \quad  \qquad \qquad
\qquad \ \ \quad (5.41)
\end{eqnarray*} where $\tilde{v}_\xi=\partial_\xi \tilde{v}$.
By Lemma 4.3 and integration by parts, we have
\begin{eqnarray*}
\displaystyle\int_{\Omega}(\nabla_{\tilde{v}}\tilde{v}+(\int_{\xi}^{1}\hbox{div}\tilde{v}d\xi^{'})\frac{\partial
\tilde{v}}{\partial
\xi})\cdot|\tilde{v}|\tilde{v}&=&\frac{1}{3}[\displaystyle\int_{\Omega}\nabla_{\tilde{v}}|\tilde{v}|^3
+\displaystyle\int_{S^2}(\int_{0}^{1}(\int_{\xi}^{1}\hbox{div} \tilde{v}d\xi^{'})d|\tilde{v}|^3)]\\
&=&\frac{1}{3}\displaystyle\int_{\Omega}(\nabla_{\tilde{v}}|\tilde{v}|^3+|\tilde{v}|^3\hbox{div}
\tilde{v})\\
&=&\frac{1}{3}\displaystyle\int_{\Omega}\hbox{div}(|\tilde{v}|^3\tilde{v})\\
&=&0. \qquad \qquad \qquad \quad\quad\quad\quad\quad\quad\quad\
(5.42)
\end{eqnarray*}
By Lemma 4.3, we get
$$\displaystyle\int_{\Omega}\hbox{div}(|\tilde{v}|^3\bar{v})
=\displaystyle\int_{\Omega}\nabla_{\bar{v}}|\tilde{v}|^3+|\tilde{v}|^3\hbox{div}\bar{v}=0.$$
By (2.19), we derive from above
$$\displaystyle\int_{\Omega}(\nabla_{\bar{v}}\tilde{v})\cdot|\tilde{v}|\tilde{v}
=\frac{1}{3}\displaystyle\int_{\Omega}\nabla_{\bar{v}}|\tilde{v}|^3=0.
\quad\quad\quad\quad\quad\quad\eqno(5.43)$$ Using Lemma 4.3, we
have
\begin{eqnarray*}
\displaystyle\int_{\Omega}\hbox{div}((|\tilde{v}|\tilde{v}\cdot
\bar{v})\tilde{v})&=&\displaystyle\int_{\Omega}\nabla_{\tilde{v}}(|\tilde{v}|\tilde{v}\cdot\bar{v})
+\displaystyle\int_{\Omega}|\tilde{v}|\tilde{v}\cdot\bar{v}
\hbox{div}\tilde{v}\\
&=&\displaystyle\int_{\Omega}(|\tilde{v}|\tilde{v}\cdot\nabla_{\tilde{v}}\bar{v}
+\bar{v}\cdot\nabla_{\tilde{v}}(|\tilde{v}|\tilde{v}))
+\displaystyle\int_{\Omega}|\tilde{v}|\tilde{v}\cdot\bar{v}\hbox{div}\tilde{v}\\
&=& 0.
\end{eqnarray*}
So
$$-\displaystyle\int_{\Omega}|\tilde{v}|\tilde{v}\cdot\nabla_{\tilde{v}}\bar{v}
=\displaystyle\int_{\Omega}\bar{v}\cdot\nabla_{\tilde{v}}(|\tilde{v}|\tilde{v})
+\displaystyle\int_{\Omega}|\tilde{v}|\tilde{v}\cdot\bar{v}\hbox{div}\tilde{v}.\eqno(5.44)$$
Using integration by parts, we obtain
\begin{eqnarray*}
\displaystyle\int_{\Omega}(\int_{0}^{1}(\tilde{v}\hbox{div}\tilde{v}
+\nabla_{\tilde{v}}\tilde{v})d\xi)\cdot|\tilde{v}|\tilde{v}
&=&\displaystyle\int_{\Omega}(\int_{0}^{1}\tilde{v}_{e_{\theta}}\tilde{v}d\xi)
\cdot\nabla_{e_{\theta}}(|\tilde{v}|\tilde{v})\\
&+&\displaystyle\int_{\Omega}(\int_{0}^{1}\tilde{v}_{e_{\varphi}}\tilde{v}d\xi)
\cdot\nabla_{e_{\varphi}}(|\tilde{v}|\tilde{v}).\qquad \quad\
(5.45)
\end{eqnarray*}
$(\frac{f}{R_0}k\times\tilde{v})\cdot|\tilde{v}|\tilde{v}=0$
implies
$$\displaystyle\int_{\Omega}(\frac{f}{R_0}k\times\tilde{v})\cdot|\tilde{v}|\tilde{v}=0.\eqno(5.46)$$
By Lemma 4.1, we get
\begin{eqnarray*}
&-&\displaystyle\int_{\Omega}[\int_{\xi}^{1}\frac{b
P}{p}\hbox{grad}((1+a
q)T)d\xi^{'}-\int_{0}^{1}\int_{\xi}^{1}\frac{b
P}{p}\hbox{grad}((1+a q)T)d\xi^{'}d\xi]\cdot |\tilde{v}|\tilde{v}\\
&=&\displaystyle\int_{\Omega}[\int_{\xi}^{1}\frac{b P}{p}(1+a
q)Td\xi^{'}-\int_{0}^{1}\int_{\xi}^{1}\frac{b P}{p}(1+a
q)Td\xi^{'}d\xi]\cdot \hbox{div}(|\tilde{v}|\tilde{v}). \ \ (5.47)
\end{eqnarray*} From (5.41) to (5.47), we get
\begin{eqnarray*}
& &\frac{1}{3}\frac{d
|\tilde{v}|_3^3}{dt}+\frac{1}{Re_1}\displaystyle\int_{\Omega}[(|\nabla_{e_{\theta}}\tilde{v}|^2
+|\nabla_{e_{\varphi}}\tilde{v}|^2)|\tilde{v}|+\frac{4}{9}|\nabla_{e_{\theta}}|\tilde{v}|^{\frac{3}{2}}|^2
+\frac{4}{9}|\nabla_{e_{\varphi}}|\tilde{v}|^{\frac{3}{2}}|^2+|\tilde{v}|^3]\\
&+&\frac{1}{Re_2}\displaystyle\int_{\Omega}(|\tilde{v}_\xi|^2|\tilde{v}|
+\frac{4}{9}|\partial_\xi|\tilde{v}|^{\frac{3}{2}}|^2)
=\displaystyle\int_{\Omega}(\bar{v}\cdot\nabla_{\tilde{v}}(|\tilde{v}|\tilde{v})
+|\tilde{v}|\tilde{v}\cdot\bar{v}\hbox{div}\tilde{v})\\
&+&\displaystyle\int_{\Omega}[(\int_{0}^{1}\tilde{v}_{\theta}\tilde{v}d\xi)\cdot
\nabla_{e_{\theta}}(|\tilde{v}|\tilde{v})
+(\int_{0}^{1}\tilde{v}_{\varphi}\tilde{v}d\xi)\cdot\nabla_{e_{\varphi}}(|\tilde{v}|\tilde{v})]\\
&+&\displaystyle\int_{\Omega}[\int_{\xi}^{1}\frac{b P}{p}(1+a
q)Td\xi^{'}-\int_{0}^{1}\int_{\xi}^{1}\frac{b P}{p}(1+a
q)Td\xi^{'}d\xi]\hbox{div}(|\tilde{v}|\tilde{v}).\qquad(5.48)
\end{eqnarray*}
By H\"older inequality, we derive from (5.48)
\begin{eqnarray*}
& &\frac{1}{3}\frac{d
|\tilde{v}|_3^3}{dt}+\frac{1}{Re_1}\displaystyle\int_{\Omega}[(|\nabla_{e_{\theta}}\tilde{v}|^2
+|\nabla_{e_{\varphi}}\tilde{v}|^2)|\tilde{v}|+\frac{4}{9}|\nabla_{e_{\theta}}|\tilde{v}|^{\frac{3}{2}}|^2
+\frac{4}{9}|\nabla_{e_{\varphi}}|\tilde{v}|^{\frac{3}{2}}|^2+|\tilde{v}|^3]\\
&+&\frac{1}{Re_2}\displaystyle\int_{\Omega}(|\tilde{v}_\xi|^2|\tilde{v}|
+\frac{4}{9}|\partial_\xi|\tilde{v}|^{\frac{3}{2}}|^2)\\
&\leq&c\displaystyle\int_{S^2}|\bar{v}|\int_{0}^{1}|\tilde{v}|^2(|\nabla_{e_{\theta}}\tilde{v}|^2
+|\nabla_{e_{\varphi}}\tilde{v}|^2)^{\frac{1}{2}}d\xi\\
&+&c\displaystyle\int_{S^2}(\int_{0}^{1}|\tilde{v}|^2d\xi)
(\int_{0}^{1}|\tilde{v}|(|\nabla_{e_{\theta}}\tilde{v}|^2
+|\nabla_{e_{\varphi}}\tilde{v}|^2)^{\frac{1}{2}}d\xi)\\
&+&c\displaystyle\int_{S^2}\overline{|T|}(\int_{0}^{1}|\tilde{v}|
(|\nabla_{e_{\theta}}\tilde{v}|^2+|\nabla_{e_{\varphi}}\tilde{v}|^2)^{\frac{1}{2}}d\xi)\\
&+& c\displaystyle\int_{S^2}\overline{|q
T|}(\int_{0}^{1}|\tilde{v}|(|\nabla_{e_{\theta}}\tilde{v}|^2
+|\nabla_{e_{\varphi}}\tilde{v}|^2)^{\frac{1}{2}}d\xi)\\
&\leq&c\displaystyle\int_{S^2}|\bar{v}|(\int_{0}^{1}|\tilde{v}|(|\nabla_{e_{\theta}}\tilde{v}|^2
+|\nabla_{e_{\varphi}}\tilde{v}|^2)d\xi)^{\frac{1}{2}}(\int_{0}^{1}|\tilde{v}|^3d\xi)^{\frac{1}{2}}\\
&+&c\displaystyle\int_{S^2}(\int_{0}^{1}|\tilde{v}|^2d\xi)(\int_{0}^{1}|\tilde{v}|d\xi)^{\frac{1}{2}}(\int_{0}^{1}
|\tilde{v}|(|\nabla_{e_{\theta}}\tilde{v}|^2+|\nabla_{e_{\varphi}}\tilde{v}|^2)d\xi)^{\frac{1}{2}}\\
&+&c\displaystyle\int_{S^2}\overline{|T|}(\int_{0}^{1}|\tilde{v}|d\xi)^{\frac{1}{2}}(\int_{0}^{1}
|\tilde{v}|(|\nabla_{e_{\theta}}\tilde{v}|^2+|\nabla_{e_{\varphi}}\tilde{v}|^2)d\xi)^{\frac{1}{2}}\\
&+&c\displaystyle\int_{S^2}\overline{|q
T|}(\int_{0}^{1}|\tilde{v}|d\xi)^{\frac{1}{2}}(\int_{0}^{1}|\tilde{v}|(|\nabla_{e_{\theta}}\tilde{v}|^2
+|\nabla_{e_{\varphi}}\tilde{v}|^2)d\xi)^{\frac{1}{2}}\\
&\leq&c\|\bar{v}\|_{L^4(S^2)}(\displaystyle\int_{\Omega}|\tilde{v}|(|\nabla_{e_{\theta}}\tilde{v}|^2
+|\nabla_{e_{\varphi}}\tilde{v}|^2))^{\frac{1}{2}}(\displaystyle\int_{S^2}
(\int_{0}^{1}|\tilde{v}|^3d\xi)^2)^{\frac{1}{4}}\\
&+&c(\displaystyle\int_{\Omega}|\tilde{v}|(|\nabla_{e_{\theta}}\tilde{v}|^2
+|\nabla_{e_{\varphi}}\tilde{v}|^2))^{\frac{1}{2}}\cdot(\displaystyle\int_{S^2}
(\int_{0}^{1}|\tilde{v}|^2d\xi)^{\frac{5}{2}})^{\frac{1}{2}}\\
&+&c\|\overline{|T|}\|_{L^4(S^2)}|\tilde{v}|_2^{\frac{1}{2}}
(\displaystyle\int_{\Omega}|\tilde{v}|(|\nabla_{e_{\theta}}\tilde{v}|^2
+|\nabla_{e_{\varphi}}\tilde{v}|^2))^{\frac{1}{2}}\\
&+&c\|\overline{|q
T|}\|_{L^4(S^2)}|\tilde{v}|_2^{\frac{1}{2}}(\displaystyle\int_{\Omega}|\tilde{v}|(|\nabla_{e_{\theta}}\tilde{v}|^2
+|\nabla_{e_{\varphi}}\tilde{v}|^2))^{\frac{1}{2}}.
\qquad\qquad\qquad\qquad\quad (5.49)
\end{eqnarray*}
By Minkowski inequality, H\"older inequality and (4.13), we have
\begin{eqnarray*}
(\displaystyle\int_{S^2}(\int_{0}^{1}|\tilde{v}|^3d\xi)^2)^\frac{1}{2}
&\leq&\int_{0}^{1}(\displaystyle\int_{S^2}(|\tilde{v}|^{\frac{3}{2}})^4)^{\frac{1}{2}}d\xi\\
&\leq&\int_{0}^{1}\||\tilde{v}|^{\frac{3}{2}}\|_{L^2(S^2)}(\|\nabla
|\tilde{v}|^{\frac{3}{2}}\|_{L^2(S^2)}^2
+\||\tilde{v}|^{\frac{3}{2}}\|_{L^2(S^2)}^2)^{\frac{1}{2}}d\xi\\
&\leq&c|\tilde{v}|_3^{\frac{3}{2}}(\int_{0}^{1}(\|\nabla
|\tilde{v}|^{\frac{3}{2}}\|_{L^2(S^2)}^2
+\||\tilde{v}|^{\frac{3}{2}}\|_{L^2(S^2)}^2)d\xi)^{\frac{1}{2}}.
\quad \  (5.50)
\end{eqnarray*}
By Minkowski inequality, H\"older inequality and
$$\|u\|_{L^5(S^2)}\leq c\|u\|^{\frac{3}{5}}_{L^3(S^2)}
\|u\|^{\frac{2}{5}}_{H^1(S^2)},\ \hbox{for any}\ u\in H^1(S^2),$$
we get
\begin{eqnarray*}
\displaystyle\int_{S^2}(\int_{0}^{1}|\tilde{v}|^2d\xi)^{\frac{5}{2}}&\leq&(\int_{0}^{1}
(\displaystyle\int_{S^2}|\tilde{v}|^5)^{\frac{2}{5}}d\xi)^{\frac{5}{2}}\\
&\leq&(\int_{0}^{1}\|\tilde{v}\|_{L^3(S^2)}^{\frac{6}{5}}
\|\tilde{v}\|_{H^1(S^2)}^{\frac{4}{5}}d\xi)^{\frac{5}{2}}\\
&\leq&c\|\tilde{v}\|^2|\tilde{v}|_3^3.\qquad \qquad \qquad
\qquad\qquad\qquad\quad\quad\ \quad\quad\quad (5.51)
\end{eqnarray*}
By (4.13), we have
$$\|\bar{v}\|_{L^4(S^2)}\leq \|\bar{v}\|_{L^2(S^2)}^{\frac{1}{2}}
\|\bar{v}\|_{H^1(S^2)}^{\frac{1}{2}}.\eqno(5.52)$$ By Minkowski
inequality, (4.13), (4.15) and H\"older inequality, we get
$$\|\overline{|T|}\|_{L^4(S^2)}=(\displaystyle\int_{S^2}(\int_{0}^{1}|T|d\xi)^4)^{\frac{1}{4}}\leq
|T|_4,\eqno{(5.53)}$$

\begin{eqnarray*}
\|\overline{|q T|}\|_{L^4(S^2)}&=&(\displaystyle\int_{S^2}(\int_{0}^{1}|q T|d\xi)^4)^{\frac{1}{4}}\\
&\leq&(\displaystyle\int_{S^2}(\int_{0}^{1}|q |^2d\xi)^2(\int_{0}^{1}|T|^2d\xi)^2)^{\frac{1}{4}}\\
&\leq&(\displaystyle\int_{S^2}(\int_{0}^{1}|q
|^2d\xi)^4)^{\frac{1}{8}}(\displaystyle\int_{S^2}(\int_{0}^{1}|T
|^2d\xi)^4)^{\frac{1}{8}}\\
&\leq&(\displaystyle\int_{0}^{1}(\int_{S^2}|q
|^8)^{\frac{1}{4}}d\xi)^{\frac{1}{2}}(\displaystyle\int_{0}^{1}(\int_{S^2}|T
|^8)^{\frac{1}{4}}d\xi)^{\frac{1}{2}}\\
&\leq&c(\displaystyle\int_{0}^{1}\|q\|_{L^4(S^2)}\|q\|_{H^1(S^2)}d\xi)^{\frac{1}{2}}
(\displaystyle\int_{0}^{1}\|T\|_{L^4(S^2)}\|T\|_{H^1(S^2)}d\xi)^{\frac{1}{2}}\\
&\leq&c|q|_4^{\frac{1}{2}}\|q\|^{\frac{1}{2}}
|T|_4^{\frac{1}{2}}\|T\|^{\frac{1}{2}}.\qquad\qquad \qquad
\qquad\qquad \qquad\quad\quad(5.54)
\end{eqnarray*}
By Young inequality, we obtain from (5.49)-(5.54)
\begin{eqnarray*}
&&\frac{d
|\tilde{v}|_3^3}{dt}+\frac{1}{Re_1}\displaystyle\int_{\Omega}[(|\nabla_{e_{\theta}}\tilde{v}|^2
+|\nabla_{e_{\varphi}}\tilde{v}|^2)|\tilde{v}|+\frac{4}{9}|\nabla_{e_{\theta}}|\tilde{v}|^{\frac{3}{2}}|^2
+\frac{4}{9}|\nabla_{e_{\varphi}}|\tilde{v}|^{\frac{3}{2}}|^2+|\tilde{v}|^3]\\
&+&\frac{1}{Re_2}\displaystyle\int_{\Omega}(|\tilde{v}_\xi|^2|\tilde{v}|
+\frac{4}{9}|\partial_\xi|\tilde{v}|^{\frac{3}{2}}|^2)\\
&\leq&c(\|\bar{v}\|_{L^2(S^2)}^2\|
\bar{v}\|_{H^1(S^2)}^2+\|\tilde{v}\|^2)|\tilde{v}|_3^3 +c(|T|_4^2
+|q|_4\|q\||T|_4\|T\|)|\tilde{v}|_2\\
&\leq&c(\|\bar{v}\|_{L^2(S^2)}^2\|
\bar{v}\|_{H^1(S^2)}^2+\|\tilde{v}\|^2)|\tilde{v}|_3^3 +c|T|_4^4
+c|T|_4^2\|T\|^2\\
&+&c(1+|q|_4^2\|q\|^2)|\tilde{v}|_2^2. \quad\quad
\quad\quad\quad\quad\quad\quad\quad\quad\quad\quad\quad\quad\quad\quad\quad\quad\quad\quad
(5.55)
\end{eqnarray*}
By Lemma 4.6, (5.17), (5.18), (5.19), (5.25), (5.39) and
$|\tilde{v}|_3^3\leq
|\tilde{v}|_2^{\frac{3}{2}}\|\tilde{v}\|_2^{\frac{3}{2}}$, we
obtain
$$|\tilde{v}(t+3r)|_3^3\leq
c(E_4+E_4^{\frac{1}{2}}E_1+2E_0(1+E_2^{\frac{1}{2}}E_1)+\frac{(4E_0E_1)^{\frac{3}{4}}}{r})\exp
c(E_0E_1+2E_1)\leq E_5,\eqno{(5.56)}$$ where $E_5$ is a positive
constant and $t\geq 0$.

\noindent{\bf $L^4$ estimates about $\tilde{v}$} \quad We take the
inner product of equation (2.23) with $|\tilde{v}|^2\tilde{v}$ in
$L^2(\Omega)\times L^2(\Omega)$, and obtain
\begin{eqnarray*}
& &\frac{1}{4}\frac{d
|\tilde{v}|_4^4}{dt}+\frac{1}{Re_1}\displaystyle\int_{\Omega}[(|\nabla_{e_{\theta}}\tilde{v}|^2
+|\nabla_{e_{\varphi}}\tilde{v}|^2)|\tilde{v}|^2+\frac{1}{2}|\nabla_{e_{\theta}}|\tilde{v}|^2|^2
+\frac{1}{2}|\nabla_{e_{\varphi}}|\tilde{v}|^2|^2+|\tilde{v}|^4]\\
&+&\frac{1}{Re_2}\displaystyle\int_{\Omega}(|\tilde{v}_\xi|^2|\tilde{v}|^2+\frac{1}{2}|\partial_\xi|\tilde{v}|^2|^2)\\
&=&-\displaystyle\int_{\Omega}(\nabla_{\tilde{v}}\tilde{v}+(\int_{\xi}^{1}\hbox{div}\tilde{v}d\xi^{'})\frac{\partial
\tilde{v}}{\partial \xi})\cdot
|\tilde{v}|^2\tilde{v}-\displaystyle\int_{\Omega}(\nabla_{\bar{v}}\tilde{v})\cdot|\tilde{v}|^2\tilde{v}
-\displaystyle\int_{\Omega}|\tilde{v}|^2\tilde{v}\cdot\nabla_{\tilde{v}}\bar{v}\\
&-&\displaystyle\int_{\Omega}[\int_{\xi}^{1}\frac{bP}{p}\hbox{grad}((1+aq)T)d\xi^{'}
-\int_{0}^{1}(\int_{\xi}^{1}\frac{bP}{p}\hbox{grad}((1+aq)T)d\xi^{'})d\xi]\cdot|\tilde{v}|^2\tilde{v}\\
&+& \displaystyle\int_{\Omega}\overline{(\tilde{v}\hbox{div}
\tilde{v}+\nabla_{\tilde{v}}\tilde{v})}\cdot|\tilde{v}|^2\tilde{v}-\displaystyle\int_{\Omega}(\frac{f}{R_0}k
\times\tilde{v})\cdot|\tilde{v}|^2\tilde{v}. \quad  \qquad \qquad
\qquad \ \ \quad (5.57)
\end{eqnarray*}
Similarly to (5.48), we derive from (5.57)
\begin{eqnarray*}
& &\frac{1}{4}\frac{d
|\tilde{v}|_4^4}{dt}+\frac{1}{Re_1}\displaystyle\int_{\Omega}[(|\nabla_{e_{\theta}}\tilde{v}|^2
+|\nabla_{e_{\varphi}}\tilde{v}|^2)|\tilde{v}|^2+\frac{1}{2}|\nabla_{e_{\theta}}|\tilde{v}|^2|^2
+\frac{1}{2}|\nabla_{e_{\varphi}}|\tilde{v}|^2|^2+|\tilde{v}|^4]\\
&+&\frac{1}{Re_2}\displaystyle\int_{\Omega}(|\tilde{v}_\xi|^2|\tilde{v}|^2+\frac{1}{2}|\partial_\xi|\tilde{v}|^2|^2)
=\displaystyle\int_{\Omega}(\bar{v}\cdot\nabla_{\tilde{v}}(|\tilde{v}|^2\tilde{v})
+|\tilde{v}|^2\tilde{v}\cdot\bar{v}\hbox{div}\tilde{v})\\
&+&\displaystyle\int_{\Omega}[(\int_{0}^{1}\tilde{v}_{\theta}\tilde{v}d\xi)\cdot
\nabla_{e_{\theta}}(|\tilde{v}|^2\tilde{v})
+(\int_{0}^{1}\tilde{v}_{\varphi}\tilde{v}d\xi)\cdot\nabla_{e_{\varphi}}(|\tilde{v}|^2\tilde{v})]\\
&+&\displaystyle\int_{\Omega}[\int_{\xi}^{1}\frac{b P}{p}(1+a
q)Td\xi^{'}-\int_{0}^{1}\int_{\xi}^{1}\frac{b P}{p}(1+a
q)Td\xi^{'}d\xi]\hbox{div}(|\tilde{v}|^2\tilde{v}).\qquad(5.58)
\end{eqnarray*}
By H\"older inequality, we derive from (5.58)
\begin{eqnarray*}
& &\frac{1}{4}\frac{d
|\tilde{v}|_4^4}{dt}+\frac{1}{Re_1}\displaystyle\int_{\Omega}[(|\nabla_{e_{\theta}}\tilde{v}|^2
+|\nabla_{e_{\varphi}}\tilde{v}|^2)|\tilde{v}|^2+\frac{1}{2}|\nabla_{e_{\theta}}|\tilde{v}|^2|^2
+\frac{1}{2}|\nabla_{e_{\varphi}}|\tilde{v}|^2|^2+|\tilde{v}|^4]\\
&+&\frac{1}{Re_2}\displaystyle\int_{\Omega}(|\tilde{v}_\xi|^2|\tilde{v}|^2+|\partial_\xi|\tilde{v}|^2|^2)\\
&\leq&c\displaystyle\int_{S^2}|\bar{v}|\int_{0}^{1}|\tilde{v}|^3(|\nabla_{e_{\theta}}\tilde{v}|^2
+|\nabla_{e_{\varphi}}\tilde{v}|^2)^{\frac{1}{2}}d\xi\\
&+&c\displaystyle\int_{S^2}(\int_{0}^{1}|\tilde{v}|^2d\xi)
(\int_{0}^{1}|\tilde{v}|^2(|\nabla_{e_{\theta}}\tilde{v}|^2
+|\nabla_{e_{\varphi}}\tilde{v}|^2)^{\frac{1}{2}}d\xi)\\
&+&c\displaystyle\int_{S^2}\overline{|T|}(\int_{0}^{1}|\tilde{v}|^2
(|\nabla_{e_{\theta}}\tilde{v}|^2+|\nabla_{e_{\varphi}}\tilde{v}|^2)^{\frac{1}{2}}d\xi)\\
&+& c\displaystyle\int_{S^2}\overline{|q
T|}(\int_{0}^{1}|\tilde{v}|^2(|\nabla_{e_{\theta}}\tilde{v}|^2
+|\nabla_{e_{\varphi}}\tilde{v}|^2)^{\frac{1}{2}}d\xi)\\
&\leq&c\displaystyle\int_{S^2}|\bar{v}|(\int_{0}^{1}|\tilde{v}|^2(|\nabla_{e_{\theta}}\tilde{v}|^2
+|\nabla_{e_{\varphi}}\tilde{v}|^2)d\xi)^{\frac{1}{2}}(\int_{0}^{1}|\tilde{v}|^4d\xi)^{\frac{1}{2}}\\
&+&c\displaystyle\int_{S^2}(\int_{0}^{1}|\tilde{v}|^2d\xi)^{\frac{3}{2}}(\int_{0}^{1}
|\tilde{v}|^2(|\nabla_{e_{\theta}}\tilde{v}|^2+|\nabla_{e_{\varphi}}\tilde{v}|^2)d\xi)^{\frac{1}{2}}\\
&+&c\displaystyle\int_{S^2}\overline{|T|}(\int_{0}^{1}|\tilde{v}|^2d\xi)^{\frac{1}{2}}(\int_{0}^{1}
|\tilde{v}|^2(|\nabla_{e_{\theta}}\tilde{v}|^2+|\nabla_{e_{\varphi}}\tilde{v}|^2)d\xi)^{\frac{1}{2}}\\
&+&c\displaystyle\int_{S^2}\overline{|q
T|}(\int_{0}^{1}|\tilde{v}|^2d\xi)^{\frac{1}{2}}(\int_{0}^{1}|\tilde{v}|^2(|\nabla_{e_{\theta}}\tilde{v}|^2
+|\nabla_{e_{\varphi}}\tilde{v}|^2)d\xi)^{\frac{1}{2}}\\
&\leq&c\|\bar{v}\|_{L^4(S^2)}(\displaystyle\int_{\Omega}|\tilde{v}|^2(|\nabla_{e_{\theta}}\tilde{v}|^2
+|\nabla_{e_{\varphi}}\tilde{v}|^2))^{\frac{1}{2}}(\displaystyle\int_{S^2}
(\int_{0}^{1}|\tilde{v}|^4d\xi)^2)^{\frac{1}{4}}\\
&+&c(\displaystyle\int_{\Omega}|\tilde{v}|^2(|\nabla_{e_{\theta}}\tilde{v}|^2
+|\nabla_{e_{\varphi}}\tilde{v}|^2))^{\frac{1}{2}}(\displaystyle\int_{S^2}
(\int_{0}^{1}|\tilde{v}|^2d\xi)^3)^{\frac{1}{2}}\\
&+&c\|\overline{|T|}\|_{L^4(S^2)}(\displaystyle\int_{\Omega}|\tilde{v}|^2(|\nabla_{e_{\theta}}\tilde{v}|^2
+|\nabla_{e_{\varphi}}\tilde{v}|^2))^{\frac{1}{2}}
(\displaystyle\int_{S^2}(\int_{0}^{1}|\tilde{v}|^2d\xi)^2)^{\frac{1}{4}}\\
&+&c\|\overline{|q
T|}\|_{L^4(S^2)}(\displaystyle\int_{\Omega}|\tilde{v}|^2(|\nabla_{e_{\theta}}\tilde{v}|^2
+|\nabla_{e_{\varphi}}\tilde{v}|^2))^{\frac{1}{2}}(\int_{S^2}(\int_{0}^{1}|\tilde{v}|^2d\xi)^2)^{\frac{1}{4}}.
\quad\quad(5.59)
\end{eqnarray*}
By Minkowski inequality, H\"older inequality and (4.13), we have
\begin{eqnarray*}
(\displaystyle\int_{S^2}(\int_{0}^{1}|\tilde{v}|^4d\xi)^2)^\frac{1}{2}
&\leq&\int_{0}^{1}(\displaystyle\int_{S^2}(|\tilde{v}|^2)^4)^{\frac{1}{2}}d\xi\\
&\leq&\int_{0}^{1}\||\tilde{v}|^2\|_{L^2(S^2)}(\|\nabla
|\tilde{v}|^2\|_{L^2(S^2)}^2
+\||\tilde{v}|^2\|_{L^2(S^2)}^2)^{\frac{1}{2}}d\xi\\
&\leq&c|\tilde{v}|_4^2(\int_{0}^{1}(\|\nabla
|\tilde{v}|^2\|_{L^2(S^2)}^2
+\||\tilde{v}|^2\|_{L^2(S^2)}^2)d\xi)^{\frac{1}{2}} \quad \quad
(5.60)
\end{eqnarray*}
By Minkowski inequality, H\"older inequality and (4.14), we get
\begin{eqnarray*}
\displaystyle\int_{S^2}(\int_{0}^{1}|\tilde{v}|^2d\xi)^3&\leq&(\int_{0}^{1}
(\displaystyle\int_{S^2}|\tilde{v}|^6)^{\frac{1}{3}}d\xi)^3\\
&\leq&(\int_{0}^{1}(\|\tilde{v}\|_{L^4(S^2)}^4
\|\tilde{v}\|_{H^1(S^2)}^2)^{\frac{1}{3}}d\xi)^3\\
&\leq&c\|\tilde{v}\|^2|\tilde{v}|_4^4.\qquad \qquad \qquad
\qquad\qquad\qquad\quad\quad\ \quad\quad(5.61)
\end{eqnarray*}
By (4.13), we have
$$\|\bar{v}\|_{L^4(S^2)}\leq \|\bar{v}\|_{L^2(S^2)}^{\frac{1}{2}}
\|\bar{v}\|_{H^1(S^2)}^{\frac{1}{2}}.\eqno(5.62)$$ By Young
inequality, (5.53) and (5.54), we obtain from (5.59)-(4.62)
\begin{eqnarray*}
& &\frac{d
|\tilde{v}|_4^4}{dt}+\frac{1}{Re_1}\displaystyle\int_{\Omega}[(|\nabla_{e_{\theta}}\tilde{v}|^2
+|\nabla_{e_{\varphi}}\tilde{v}|^2)|\tilde{v}|^2+\frac{1}{2}|\nabla_{e_{\theta}}|\tilde{v}|^2|^2
+\frac{1}{2}|\nabla_{e_{\varphi}}|\tilde{v}|^2|^2+|\tilde{v}|^4]\\
&+&\frac{1}{Re_2}\displaystyle\int_{\Omega}(|\tilde{v}_\xi|^2|\tilde{v}|^2+\frac{1}{2}|\partial_\xi|\tilde{v}|^2|^2)\\
&\leq&c\|\bar{v}\|_{L^2(S^2)}^2\|
\bar{v}\|_{H^1(S^2)}^2|\tilde{v}|_4^4
+c\|\tilde{v}\|^2|\tilde{v}|_4^4 +c|T|_4^2 |\tilde{v}|_4^2
+c|q|_4\|q\||T|_4\|T\||\tilde{v}|_4^2\\
&\leq&c(\|\bar{v}\|_{L^2(S^2)}^2\|\bar{v}\|_{H^1(S^2)}^2+|T|_4^2
+\|\tilde{v}\|^2 +|q|_4^2\|q\|^2)|\tilde{v}|_4^4\\
&+&c|T|_4^2+c|T|_4^2\|T\|^2.
 \quad\quad\quad\quad\quad\quad\quad\quad\quad\quad\quad\quad\quad\quad
 \quad\quad\quad\quad\quad\quad(5.63)
\end{eqnarray*}
By Lemma 4.6, (5.17), (5.18), (5.19), (5.25), (5.39), (5.56) and
$|\tilde{v}|_4^4\leq |\tilde{v}|_3^2\|\tilde{v}\|^2$, we obtain
$$|\tilde{v}(t+4r)|_4^4\leq
c(E_4^{\frac{1}{2}}+E_1E_4^{\frac{1}{2}}+\frac{E_5^{\frac{2}{3}}E_1}{r})\exp
c(E_0E_1+E_1+E_4^{\frac{1}{2}}+E_2^{\frac{1}{2}}E_1)\leq
E_6,\eqno{(5.64)}$$ where $E_6$ is a positive constant and $t\geq
0$. From (5.63) and (5.64), we have
\begin{eqnarray*}
& &\int_{t+4r}^{t+5r}[\frac{1}{Re_1}\displaystyle\int_{\Omega}
((|\nabla_{e_{\theta}}\tilde{v}|^2
+|\nabla_{e_{\varphi}}\tilde{v}|^2)|\tilde{v}|^2+\frac{1}{2}|\nabla_{e_{\theta}}|\tilde{v}|^2|^2
+\frac{1}{2}|\nabla_{e_{\varphi}}|\tilde{v}|^2|^2+|\tilde{v}|^4)\\
&+&\frac{1}{Re_2}\displaystyle\int_{\Omega}(|\tilde{v}_\xi|^2|\tilde{v}|^2
+\frac{1}{2}|\partial_\xi|\tilde{v}|^2|^2)]\leq
E_6^2+E_6=E_7. \quad\quad\quad\quad\quad\quad\quad\quad (5.65)
\end{eqnarray*}
By Gronwall inequality, from (5.63) we obtain
$$|\tilde{v}(t)|_4^4\leq
C_2,\eqno{(5.66)}$$ where $C_2=C_2(\|U_0\|,\ \|Q_1\|_1,\
\|Q_2\|_1)>0$ and $0\leq t<4r $.

\noindent{\bf $H^1$ estimates about $\bar{v}$} \quad Taking the
inner product of the equation (2.22) with $-\triangle \bar{v}$ in
$L^2(S^2)\times L^2(S^2)$, we get
$$\frac{1}{2}\frac{d
\|\bar{v}\|_{H^1(S^2)}^2}{dt}+\frac{1}{Re_1}\|\triangle
\bar{v}\|_{L^2(S^2)}^2=\displaystyle\int_{S^2}[\nabla_{\bar{v}}\bar{v}+\int_{0}^{1}(\tilde{v}\hbox{div}\tilde{v}
+\nabla_{\tilde{v}}\tilde{v})d\xi]\cdot\triangle\bar{v}$$
$$+\displaystyle\int_{S^2}(\hbox{grad}\Phi_{s}+\frac{f}{R_0}k\times
\bar{v})\cdot\triangle\bar{v}+\displaystyle\int_{S^2}[\int_{0}^{1}\int_{\xi}^{1}\frac{b
P}{p}\hbox{grad}((1+aq)T)d\xi^{'}d\xi]\cdot\triangle\bar{v}.
\eqno{(5.67)}
$$
By H\"older inequality, (4.13) and Young inequality, we have
\begin{eqnarray*}
&
&|\displaystyle\int_{S^2}(\nabla_{\bar{v}}\bar{v}\cdot\triangle\bar{v})|\\
&\leq&c\displaystyle\int_{S^2}|\bar{v}|(|\nabla_{e_{\theta}}\bar{v}|^2
+|\nabla_{e_{\varphi}}\bar{v}|^2)^{\frac{1}{2}}|\triangle\bar{v}|\\
&\leq&c\|\bar{v}\|_{L^4(S^2)}(\displaystyle\int_{S^2}(|\nabla_{e_{\theta}}\bar{v}|^2
+|\nabla_{e_{\varphi}}\bar{v}|^2)^2)^{\frac{1}{4}}\|\triangle\bar{v}\|_{L^2(S^2)}\\
&\leq&c\|\bar{v}\|_{L^2(S^2)}^{\frac{1}{2}}\|\bar{v}\|_{H^1(S^2)}^{\frac{1}{2}}
(\displaystyle\int_{S^2}(|\nabla_{e_{\theta}}\bar{v}|^2+|\nabla_{e_{\varphi}}\bar{v}|^2))^{\frac{1}{4}}\\
&\times&(\displaystyle\int_{S^2}(|\nabla_{e_{\theta}}\bar{v}|^2
+|\nabla_{e_{\varphi}}\bar{v}|^2)+\|\triangle\bar{v}\|_{L^2(S^2)}^2)^{\frac{1}{4}}\|\triangle\bar{v}\|_{L^2(S^2)}\\
&\leq&c\|\bar{v}\|_{L^2(S^2)}^{\frac{1}{2}}\|\bar{v}\|_{H^1(S^2)}^{\frac{1}{2}}
(\displaystyle\int_{S^2}(|\nabla_{e_{\theta}}\bar{v}|^2+|\nabla_{e_{\varphi}}\bar{v}|^2))^{\frac{1}{4}}\\
&\times&[(\displaystyle\int_{S^2}(|\nabla_{e_{\theta}}\bar{v}|^2
+|\nabla_{e_{\varphi}}\bar{v}|^2))^{\frac{1}{4}}+
\|\triangle\bar{v}\|_{L^2(S^2)}^{\frac{1}{2}}]\|\triangle\bar{v}\|_{L^2(S^2)}\\
&\leq&c\|\bar{v}\|_{L^2(S^2)}^{\frac{1}{2}}\|\bar{v}\|_{H^1(S^2)}^{\frac{3}{2}}
\|\triangle\bar{v}\|_{L^2(S^2)}+c\|\bar{v}\|_{L^2(S^2)}^{\frac{1}{2}}
\|\bar{v}\|_{H^1(S^2)}\|\triangle\bar{v}\|_{L^2(S^2)}^{\frac{3}{2}}\\
&\leq&c\|\bar{v}\|_{L^2(S^2)}\|\bar{v}\|_{H^1(S^2)}^3+c\|\bar{v}\|_{L^2(S^2)}^2\|\bar{v}\|_{H^1(S^2)}^4
+\varepsilon \|\triangle\bar{v}\|^2_{L^2(S^2)}\\
&\leq&c(\|\bar{v}\|_{L^2(S^2)}^2+\|\bar{v}\|_{H^1(S^2)}^2
+\|\bar{v}\|_{L^2(S^2)}^2\|\bar{v}\|_{H^1(S^2)}^2)\|\bar{v}\|_{H^1(S^2)}^2+\varepsilon
\|\triangle\bar{v}\|^2_{L^2(S^2)}.(5.68)
\end{eqnarray*}
By H\"older inequality and Minkowski inequality, we obtain
\begin{eqnarray*}
&
&|\displaystyle\int_{S^2}(\int_{0}^{1}(\tilde{v}\hbox{div}\tilde{v}
+\nabla_{\tilde{v}}\tilde{v})d\xi\cdot\triangle \bar{v})|\\
&\leq&\displaystyle\int_{S^2}\int_{0}^{1}|\tilde{v}|(|\nabla_{e_{\theta}}\tilde{v}|^2
+|\nabla_{e_{\varphi}}\tilde{v}|^2)^\frac{1}{2}d\xi|\triangle\bar{v}|\\
&\leq&[\displaystyle\int_{S^2}(\int_{0}^{1}|\tilde{v}|(|\nabla_{e_{\theta}}\tilde{v}|^2
+|\nabla_{e_{\varphi}}\tilde{v}|^2)^\frac{1}{2}d\xi)^2]^{\frac{1}{2}}\|\triangle\bar{v}\|_{L^2(S^2)}\\
&\leq&\displaystyle\int_{0}^{1}(\displaystyle\int_{S^2}|\tilde{v}|^2(|\nabla_{e_{\theta}}\tilde{v}|^2
+|\nabla_{e_{\varphi}}\tilde{v}|^2))^\frac{1}{2}d\xi\|\triangle\bar{v}\|_{L^2(S^2)}\\
&\leq&c\displaystyle\int_{\Omega}|\tilde{v}|^2(|\nabla_{e_{\theta}}\tilde{v}|^2
+|\nabla_{e_{\varphi}}\tilde{v}|^2)+\varepsilon\|\triangle\bar{v}\|_{L^2(S^2)}^2.\qquad
\quad\qquad \qquad\qquad\ \ \ \ (5.69)
\end{eqnarray*}
Using Lemma 4.1 and (2.19), we have
$$\displaystyle\int_{S^2}\hbox{grad}\Phi_s\cdot\triangle\bar{v}=0,$$
$$\displaystyle\int_{S^2}[\int_{0}^{1}\int_{\xi}^{1}\frac{b
P}{p}\hbox{grad}((1+aq)T)d\xi^{'}d\xi]\cdot\triangle\bar{v}=0.\qquad
\qquad\qquad\eqno(5.70)$$ $(\frac{f}{R_0}k\times
\bar{v})\cdot\triangle\bar{v}=0$ implies
$$\displaystyle\int_{S^2}(\frac{f}{R_0}k\times
\bar{v})\cdot\triangle\bar{v}=0.\eqno(5.71)$$ From (5.67)-(5.71),
choosing $\varepsilon$ small enough, we obtain
\begin{eqnarray*}
& &\frac{d \|\bar{v}\|_{H^1(S^2)}^2}{dt}+\frac{1}{Re_1}\|\triangle
\bar{v}\|_{L^2(S^2)}^2\\&\leq&c(\|\bar{v}\|_{L^2(S^2)}^2+\|\bar{v}\|_{H^1(S^2)}^2
+\|\bar{v}\|_{L^2(S^2)}^2\|\bar{v}\|_{H^1(S^2)}^2)\|\bar{v}\|_{H^1(S^2)}^2\\
&+&c\displaystyle\int_{\Omega}|\tilde{v}|^2(|\nabla_{e_{\theta}}\tilde{v}|^2
+|\nabla_{e_{\varphi}}\tilde{v}|^2).\qquad \qquad
\qquad\qquad\qquad\quad\quad \quad\quad\quad\quad(5.72)
\end{eqnarray*}
By Lemma 4.6, (5.17), (5.18), (5.19) and (5.65), we get
$$\|\bar{v}(t+5r)\|_{H^1(S^2)}^2\leq c(\frac{E_1}{r}+E_7)\exp
c(E_0E_1+E_1 )\leq E_8,\eqno{(5.73)}$$ where $E_8$ is a positive
constant. By Gronwall inequality, from (5.72) we obtain
$$\|\bar{v}(t)\|_{H^1(S^2)}^2\leq
C_3,\eqno{(5.74)}$$ where $C_3=C_3(\|U_0\|,\ \|Q_1\|_1,\
\|Q_2\|_1)>0$ and $0\leq t<5r $.

\smallskip
\noindent{\bf $L^2$ estimates about $v_{\xi}$} \quad Taking the
derivative, with respect to $\xi$, of equation (2.11), we get the
following equation
\begin{eqnarray*}
& &\frac{\partial v_{\xi}}{\partial t}-\frac{1}{Re_1}\triangle
v_{\xi}-\frac{1}{Re_2}\frac{\partial^2 v_{\xi}}{\partial
\xi^2}+\nabla_{v}v_{\xi}+(\displaystyle\int_{\xi}^{1}\hbox{div}
vd\xi^{'})\frac{\partial v_{\xi}}{\partial \xi}\\
&+&\nabla_{v_{\xi}}v-(\hbox{div} v)\frac{\partial v}{\partial
\xi}+\frac{f}{R_0}k\times
v_{\xi}-\frac{bP}{p}\hbox{grad}[(1+aq)T]=0. \qquad \quad\ (5.75)
\end{eqnarray*}
Taking the inner product of equation (5.75) with $v_{\xi}$ in
$L^2(\Omega)\times L^2(\Omega)$, we obtain

\begin{eqnarray*}
& &\frac{1}{2}\frac{d |v_{\xi}|_2^2}{d
t}+\frac{1}{Re_1}\displaystyle\int_{\Omega}(|\nabla_{e_{\theta}}v_{\xi}|^2
+|\nabla_{e_{\varphi}}v_{\xi}|^2+|v_{\xi}|^2)
+\frac{1}{Re_2}\displaystyle\int_{\Omega}|\frac{\partial
v_{\xi}}{\partial \xi}|^2\\
&=&-\displaystyle\int_{\Omega}(\nabla_{v}v_{\xi}+(\int_{\xi}^{1}\hbox{div}
vd\xi^{'})\frac{\partial v_{\xi}}{\partial \xi})\cdot
v_{\xi}-\displaystyle\int_{\Omega}(\nabla_{v_{\xi}}v-(\hbox{div}
v)\frac{\partial v}{\partial \xi})\cdot
v_{\xi}\\&-&\displaystyle\int_{\Omega}(\frac{f}{R_0}k\times
v_{\xi})\cdot
v_{\xi}+\displaystyle\int_{\Omega}\frac{bP}{p}\hbox{grad}[(1+aq)T]\cdot
v_{\xi}.\qquad\qquad \qquad \ \ \ (5.76)
\end{eqnarray*}
By integration by parts and Lemma 4.3, we get
$$\displaystyle\int_{\Omega}(\nabla_{v}v_{\xi}+(\int_{\xi}^{1}\hbox{div}
vd\xi^{'})\frac{\partial v_{\xi}}{\partial \xi})\cdot v_{\xi}=0.
\eqno(5.77)$$ By integration by parts, H\"older inequality, (4.16)
and Young inequality, we have
\begin{eqnarray*}
-\displaystyle\int_{\Omega}(\nabla_{v_{\xi}}v-(\hbox{div}
v)\frac{\partial v}{\partial \xi})\cdot
v_{\xi}&\leq&c\displaystyle\int_{\Omega}|v||v_{\xi}|(|\nabla_{e_{\theta}}v_{\xi}|^2
+|\nabla_{e_{\varphi}}v_{\xi}|^2)^{\frac{1}{2}}\\
&\leq&c|v|_4|v_{\xi}|_4(\displaystyle\int_{\Omega}(|\nabla_{e_{\theta}}v_{\xi}|^2
+|\nabla_{e_{\varphi}}v_{\xi}|^2))^{\frac{1}{2}}\\
&\leq&c|v|_4|v_{\xi}|_2^{\frac{1}{4}}\|v_{\xi}\|^{\frac{3}{4}}
(\displaystyle\int_{\Omega}(|\nabla_{e_{\theta}}v_{\xi}|^2+|\nabla_{e_{\varphi}}v_{\xi}|^2))^{\frac{1}{2}}\\
&\leq&\varepsilon\|v_{\xi}\|^2+c|v|_4^8|v_{\xi}|_2^2.\qquad
\qquad\qquad\qquad(5.78)
\end{eqnarray*}
$(\frac{f}{R_0}k\times v_{\xi})\cdot v_{\xi}=0$ implies
$$\displaystyle\int_{\Omega}(\frac{f}{R_0}k\times v_{\xi})\cdot v_{\xi}=0.\eqno(5.79)$$
By Lemma 4.1, H\"older inequality and Young inequality, we obtain
\begin{eqnarray*}
\displaystyle\int_{\Omega}\frac{bP}{p}\hbox{grad}[(1+aq)T]\cdot
v_{\xi}&=&-\displaystyle\int_{\Omega}\frac{bP}{p}(1+aq)T
\hbox{div}v_{\xi}\\
&\leq&c\displaystyle\int_{\Omega}|T||\hbox{div}
v_{\xi}|+c\displaystyle\int_{\Omega}|q||T||\hbox{div} v_{\xi}|\\
&\leq&c|T|_2^2+c|q|_4^2|T|_4^2+\varepsilon\|v_{\xi}\|^2.
\quad\quad \quad \quad\quad \ (5.80)
\end{eqnarray*}
Choosing $\varepsilon$ small enough, we derive from (5.76)-(5.80),
\begin{eqnarray*}
& &\frac{d |v_{\xi}|_2^2}{d
t}+\frac{1}{Re_1}\displaystyle\int_{\Omega}(|\nabla_{e_{\theta}}v_{\xi}|^2
+|\nabla_{e_{\varphi}}v_{\xi}|^2+|v_{\xi}|^2)
+\frac{1}{Re_2}\displaystyle\int_{\Omega}|\frac{\partial
v_{\xi}}{\partial \xi}|^2\\
&\leq&c|v|_4^8|v_{\xi}|_2^2+c|T|_2^2+c|q|_4^2|T|_4^2\\
&\leq&c|\bar{v}+\tilde{v}|_4^8|v_{\xi}|_2^2+c|T|_2^2+c|q|_4^2|T|_4^2\\
&\leq&c(|\bar{v}|_{H^1(S^2)}^8+|\tilde{v}|_4^8)|v_{\xi}|_2^2+c|T|_2^2+c|q|_4^4+c|T|_4^4.\qquad\qquad\qquad\quad\
\quad (5.81)
\end{eqnarray*}
By Lemma 4.6, (5.18), (5.25), (5.39), (5.64) and (5.73), we get
$$|v_{\xi}(t+6r)|_2^2\leq c(E_0+\frac{E_1}{r}+E_2+E_4)\exp c(E_6^2+E_8^4
)\leq E_9,\eqno{(5.82)}$$ where $E_9$ is a positive constant and
$t\geq 0$. From (5.81) and (5.82), we have
$$c_1\int_{t+6r}^{t+7r}\|v_{\xi}\|^2\leq E_9^2+E_9=E_{10}.\eqno{(5.83)}$$ By
Gronwall inequality, from (5.81) we obtain
$$|v_{\xi}(t)|_2^2\leq
C_4,\eqno{(5.84)}$$ where  $C_4=C_4(\|U_0\|,\ \|Q_1\|_1,\
\|Q_2\|_1)>0$ and $0\leq t<6r $.

\noindent{\bf $L^2$ estimates about $T_{\xi}$ $q_{\xi}$}\quad
Taking the derivative, with respect to $\xi$, of the equations
(2.12), (2.13), we get the following equations

\begin{eqnarray*}
& &\frac{\partial T_{\xi}}{\partial t}-\frac{1}{Rt_1}\triangle
T_{\xi}-\frac{1}{Rt_2}\frac{\partial^2 T_{\xi}}{\partial
\xi^2}+\nabla_{v}T_{\xi}+W(v)\frac{\partial T_{\xi}}{\partial
\xi}+ \nabla_{v_{\xi}}T -(\hbox{div} v)\frac{\partial T}{\partial
\xi}\\
&+&\frac{bP}{p}(1+aq)\hbox{div}
v-\frac{abP}{p}q_{\xi}W(v)+\frac{bP(P-p_0)}{p^2}(1+aq)W(v)=Q_{1\xi},(5.85)
\end{eqnarray*}
$$\frac{\partial q_{\xi}}{\partial t}-\frac{1}{Rq_1}\triangle
q_{\xi}-\frac{1}{Rq_2}\frac{\partial^2 q_{\xi}}{\partial
\xi^2}+\nabla_{v}q_{\xi}+W(v)\frac{\partial q_{\xi}}{\partial \xi}
+\nabla_{v_{\xi}}q-(\hbox{div} v)\frac{\partial q}{\partial \xi}
=Q_{2\xi}. \eqno(5.86)$$ We take the inner product of equation
(5.85) with $T_{\xi}$ in $L^2(\Omega)$ and obtain
\begin{eqnarray*}
&&\frac{1}{2}\frac{d |T_{\xi}|_2^2}{d
t}+\frac{1}{Rt_1}\displaystyle\int_{\Omega}|\nabla
T_{\xi}|^2+\frac{1}{Rt_2}\displaystyle\int_{\Omega}|
T_{\xi\xi}|^2-\frac{1}{Rt_2}\displaystyle\int_{S^2}(T_\xi|_{\xi=1}\cdot
T_{\xi\xi}|_{\xi=1})\\
&=&-\displaystyle\int_{\Omega}(\nabla_{v}T_{\xi}+W(v)\frac{\partial
T_{\xi}}{\partial
\xi})T_{\xi}-\displaystyle\int_{\Omega}(\nabla_{v_{\xi}}T-(\hbox{div}v)\frac{\partial
T}{\partial \xi})T_{\xi}+\displaystyle\int_{\Omega}Q_{1\xi}T_{\xi}\\
&+&\displaystyle\int_{\Omega}[-\frac{b P}{p}(1+a q)(\hbox{div}
v)-\frac{b P(P-p_0)}{p}(1+a q)W(v)+\frac{ab P}{p}
q_{\xi}W(v)]T_{\xi}.(5.87)
\end{eqnarray*}
Similarly to (5.77), we have
$$-\displaystyle\int_{\Omega}(\nabla_{v}T_{\xi}+W(v)\frac{\partial
T_{\xi}}{\partial \xi})T_{\xi}=0. \eqno(5.88)$$ By integration by
parts, H\"older inequality, (4.16), Poincar\'e inequality and
Young inequality, we obtain
\begin{eqnarray*}
&
&|\displaystyle\int_{\Omega}(\nabla_{v_{\xi}}T-\hbox{div}(v)\frac{\partial
T}{\partial \xi})T_{\xi}|\\
&\leq&c\displaystyle\int_{\Omega}[(|\nabla_{e_{\theta}}v_{\xi}|^2
+|\nabla_{e_{\varphi}}v_{\xi}|^2)^{\frac{1}{2}}|T||T_{\xi}|+|v_{\xi}||T||\nabla
T_{\xi}|+|v||\nabla T_{\xi}||T_{\xi}|]\\
&\leq&c(\displaystyle\int_{\Omega}(|\nabla_{e_{\theta}}v_{\xi}|^2
+|\nabla_{e_{\varphi}}v_{\xi}|^2))^{\frac{1}{2}}|T|_4|T_{\xi}|_4+c|v_{\xi}|_4|T|_4|\nabla
T_{\xi}|_2+c|v|_4|\nabla T_{\xi}|_2|T_{\xi}|_4\\
&\leq&c\displaystyle\int_{\Omega}(|\nabla_{e_{\theta}}v_{\xi}|^2
+|\nabla_{e_{\varphi}}v_{\xi}|^2)+\frac{\varepsilon}{2}|\nabla
T_{\xi}|_2^2+c|T|_4^2|T_{\xi}|_4^2+c|v_{\xi}|_4^2|T|_4^2+c|v|_4^2|T_{\xi}|_4^2\\
&\leq&c\displaystyle\int_{\Omega}(|\nabla_{e_{\theta}}v_{\xi}|^2
+|\nabla_{e_{\varphi}}v_{\xi}|^2)+\frac{\varepsilon}{2}|\nabla
T_{\xi}|_2^2+c(|T|_4^2+|v|_4^2)|T_{\xi}|_2^{\frac{1}{2}}(|\nabla
T_{\xi}|_2^2+|T_{\xi\xi}|_2^2)^{\frac{3}{4}}\\
&+&c|T|_4^2|v_{\xi}|_2^{\frac{1}{2}}(|v_{\xi\xi}|_2^2
+\displaystyle\int_{\Omega}(|\nabla_{e_{\theta}}v_{\xi}|^2+|\nabla_{e_{\varphi}}v_{\xi}|^2))^{\frac{3}{4}}\\
&\leq&\varepsilon(|T_{\xi\xi}|_2^2+|\nabla T_{\xi}|_2^2)+c(|v_{\xi
\xi}|_2^2+\displaystyle\int_{\Omega}(|\nabla_{e_{\theta}}v_{
\xi}|^2+|\nabla_{e_{\varphi}}v_{\xi}|^2))\\
&+&c|T|_4^8|v_{\xi}|_2^2+c(|T|_4^8+|v|_4^8)|T_{\xi}|_2^2.
\qquad\qquad\qquad\qquad\qquad\qquad\qquad\qquad\ (5.89)
\end{eqnarray*}
By integrating by parts, H\"older inequality, Minkowsky
inequality, Poincar\'e inequality, Young inequality and Lemma 4.2,
we obtain
\begin{eqnarray*}
& &|\displaystyle\int_{\Omega}[-\frac{b P}{p}(1+a q)(\hbox{div}
v)-\frac{b P(P-p_0)}{p}(1+a q)W(v)+\frac{ab P}{p}
q_{\xi}W(v)]T_{\xi}|\\
&\leq&|\displaystyle\int_{\Omega}(\frac{b P}{p}\nabla q\cdot
vT_{\xi}+\frac{b P}{p}(1+a q)v\cdot\nabla
T_{\xi})|+|\displaystyle\int_{\Omega}[\frac{b P(P-p_0)}{p}(1+a
q)W(v)]T_{\xi}|\\
&+&|\displaystyle\int_{\Omega}(\frac{ab P}{p}(\nabla
q_{\xi})(\int_{\xi}^{1}v)T_{\xi}+\frac{ab
P}{p}q_{\xi}(\int_{\xi}^{1}v)\nabla T_{\xi})|\\
&\leq&c(\displaystyle\int_{\Omega}|\nabla
q|^2)^{\frac{1}{2}}|v|_4|T_{\xi}|_4+c|q|_4|v|_4|\nabla
T_{\xi}|_2+c|v|_2|\nabla T_{\xi}|_2\\
&+&c(|T_{\xi}|_2\|v\|+|q|_4|T_{\xi}|_4\|v\|)+
c|\nabla q_{\xi}|_2|v|_4|T_{\xi}|_4+c|q_{\xi}|_4|v|_4|\nabla T_{\xi}|_2\\
&\leq&\varepsilon(|\nabla q_{\xi}|_2^2+|\nabla
T_{\xi}|_2^2)+c|\nabla
q|_2^2+c\|v\|^2+c(|v|_4^2+|q|_4^2)|T_{\xi}|_4^2\\
&+&c|v|_4^2(|q_{\xi}|_4^2+|q|_4^2)+c|T_{\xi}|_2^2\\
&\leq&\varepsilon(|\nabla q_{\xi}|_2^2+|\nabla
T_{\xi}|_2^2)+c|\nabla
q|_2^2+c\|v\|^2+c(|v|_4^2+|q|_4^2)|T_{\xi}|_2^{\frac{1}{2}}\|T_{\xi}\|^{\frac{3}{2}}
\\
&+&c|v|_4^2|q_{\xi}|_2^{\frac{1}{2}}\|q_{\xi}\|^{\frac{3}{2}}+c|v|_4^2|q|_4^2+c|T_{\xi}|_2^2\\
&\leq&\varepsilon(|\nabla T_{\xi}|_2^2+|T_{\xi,
\xi}|_2^2)+\varepsilon(|\nabla q_{\xi}|_2^2+|q_{\xi,
\xi}|_2^2)+c(|v|_4^8+|q|_4^8+1)|T_{\xi}|_2^2\\
&+&c|v|_4^8|q_{\xi}|_2^2+c(|\nabla
q|_2^2+\|v\|^2)+c|v|_4^2|q|_4^2. \qquad \qquad \qquad\quad
\qquad\quad \qquad(5.90)
\end{eqnarray*}
By integration by parts, H\"older inequality and Young inequality,
we have
$$|\displaystyle\int_{\Omega}Q_{1\xi}T_{\xi}|\leq c|Q_{1\xi}|_2^2+c|T_{\xi}|_2^2.\eqno(5.91)$$
Taking the trace on $\xi=1$ of equation (2.12), we have
$$\frac{1}{Rt_2}T_{\xi\xi}|_{\xi=1}=\frac{\partial T|_{\xi=1}}{\partial t}
+(\nabla_v T)|_{\xi=1}-\frac{1}{Rt_1}\triangle
T|_{\xi=1}-Q_1|_{\xi=1}. \eqno(5.92)$$ From (2.15), (5.92), we get

\begin{eqnarray*}
& &-\frac{1}{Rt_2}\displaystyle\int_{S^2}(T_\xi|_{\xi=1}
T_{\xi\xi}|_{\xi=1})\\
&=&\alpha_s\displaystyle\int_{S^2}T|_{\xi=1}(\frac{\partial
T|_{\xi=1}}{\partial t} +(\nabla_v
T)|_{\xi=1}-\frac{1}{Rt_1}\triangle T|_{\xi=1}-Q_1|_{\xi=1})\\
&=&\alpha_s(\frac{1}{2}\frac{d|T|_{\xi=1}|^2_2}{dt}+\frac{1}{Rt_1}|\nabla
T|_{\xi=1}|^2_2)+\alpha_s\displaystyle\int_{S^2}T|_{\xi=1}(
(\nabla_v T)|_{\xi=1}-Q_1|_{\xi=1}). (5.93)
\end{eqnarray*}
By Lemma 4.3, we have

\begin{eqnarray*}
& &-\alpha_s\displaystyle\int_{S^2}T|_{\xi=1}( (\nabla_v
T)|_{\xi=1}-Q_1|_{\xi=1})\\
&=&-\frac{\alpha_s}{2}\displaystyle\int_{S^2}(\nabla_v
T^2)|_{\xi=1}+\alpha_s\displaystyle\int_{S^2}T|_{\xi=1}
Q_1|_{\xi=1}\\
&=&\frac{\alpha_s}{2}\displaystyle\int_{S^2}T^2|_{\xi=1}
\hbox{div}v|_{\xi=1}+\alpha_s\displaystyle\int_{S^2}T|_{\xi=1}
Q_1|_{\xi=1}\\
&=&\frac{\alpha_s}{2}\displaystyle\int_{S^2}T^2|_{\xi=1}
(\displaystyle\int^1_\xi\hbox{div}v_{\xi}d\xi'+\hbox{div}v)+\alpha_s\displaystyle\int_{S^2}T|_{\xi=1}
Q_1|_{\xi=1}\\
&\leq&c|T|_{\xi=1}|^4_4+c\|v_{\xi}\|^2+c\|v\|^2+c|T|_{\xi=1}|^2_2+c|Q_1|_{\xi=1}|^2_2.
\quad\quad\quad\quad\quad\quad(5.94)
\end{eqnarray*}

We derive from (5.87)-(5.94)
\begin{eqnarray*}
& &\frac{1}{2}\frac{d (|T_{\xi}|_2^2+\alpha_s|T|_{\xi=1}|^2_2)}{d
t}+\frac{1}{Rt_1}\displaystyle\int_{\Omega}|\nabla
T_{\xi}|^2+\frac{1}{Rt_2}\displaystyle\int_{\Omega}|
T_{\xi\xi}|^2+\frac{\alpha_s}{Rt_1}|\nabla
T|_{\xi=1}|^2_2\\
&\leq&2\varepsilon(|T_{\xi\xi}|_2^2+|\nabla
T_{\xi}|_2^2)+\varepsilon(|\nabla q_{\xi}|_2^2+|q_{\xi\xi}|_2^2)
+c(1+|T|_4^8+|v|_4^8+|q|_4^8)|T_{\xi}|_2^2\\
&+&c|v|_4^8|q_{\xi}|_2^2+c\|v_{\xi}\|^2+c\|v\|^2+c|\nabla
q|_2^2+c|T|_4^8|v_{\xi}|_2^2+c|q|_4^2|v|_4^2\\
&+&c|T|_{\xi=1}|^4_4+c|T|_{\xi=1}|^2_2+c|Q_1|_{\xi=1}|^2_2+c|Q_1|_2^2+c|Q_{1\xi}|_2^2.
\quad\quad\quad\quad\quad\quad(5.95)
\end{eqnarray*}

Similarly to (5.95), we have
\begin{eqnarray*}
& &\frac{1}{2}\frac{d (|q_{\xi}|_2^2+\beta_s|q|_{\xi=1}|^2_2)}{d
t}+\frac{1}{Rq_1}\displaystyle\int_{\Omega}|\nabla
q_{\xi}|^2+\frac{1}{Rq_2}\displaystyle\int_{\Omega}|
q_{\xi\xi}|^2+\frac{\beta_s}{Rq_1}|\nabla
q|_{\xi=1}|^2_2\\
&\leq&\varepsilon(|\nabla q_{\xi}|_2^2+|q_{\xi\xi}|_2^2)
+c(|v|_4^8+|T|_4^8)|q_{\xi}|_2^2\\
&+&c\|v_{\xi}\|^2+c\|v\|^2+c|q|_4^8|v_{\xi}|_2^2\\
&+&c|q|_{\xi=1}|^4_4+c|q|_{\xi=1}|^2_2+c|Q_2|_{\xi=1}|^2_2+c|Q_2|_2^2+c|Q_{2\xi}|_2^2.
\quad\quad\quad\quad\quad\quad\quad (5.96)
\end{eqnarray*}
From (5.95) and (5.96), choosing $\varepsilon$ small enough, we
obtain
\begin{eqnarray*}
& &\frac{d
(|T_{\xi}|_2^2+|q_{\xi}|_2^2+\beta_s|q|_{\xi=1}|^2_2+\alpha_s|T|_{\xi=1}|^2_2)}{d
t}+\frac{1}{Rt_1}\displaystyle\int_{\Omega}|\nabla
T_{\xi}|^2+\frac{1}{Rt_2}\displaystyle\int_{\Omega}|
T_{\xi\xi}|^2\\
&+&\frac{\alpha_s}{Rt_1}|\nabla
T|_{\xi=1}|^2_2+\frac{1}{Rq_1}\displaystyle\int_{\Omega}|\nabla
q_{\xi}|^2+\frac{1}{Rq_2}\displaystyle\int_{\Omega}|
q_{\xi\xi}|^2+\frac{\beta_s}{Rq_1}|\nabla
q|_{\xi=1}|^2_2\\
&\leq&c(1+|T|_4^8+|v|_4^8+|q|_4^8)(|T_{\xi}|_2^2+|q_{\xi}|_2^2)+c\|v_{\xi}\|^2+c\|v\|^2
+c\|q\|^2\\
&+&c(|T|_4^8+|q|_4^8)|v_{\xi}|_2^2+c|q|_4^2|v|_4^2+c|T|_{\xi=1}|^4_4
+c|T|_{\xi=1}|^2_2+c|q|_{\xi=1}|^4_4+c|q|_{\xi=1}|^2_2\\
&+&c(|Q_1|_{\xi=1}|^2_2+|Q_2|_{\xi=1}|^2_2)+c(|Q_1|_2^2+|Q_{1\xi}|_2^2+|Q_2|_2^2+|Q_{2\xi}|_2^2).
\quad\quad\quad\ (5.97)
\end{eqnarray*}

By Lemma 4.6, (5.18), (5.25), (5.26), (5.39), (5.40), (5.64),
(5.73), (5.82), (5.83), we get
$$|T_{\xi}(t+7r)|_2^2+|q_{\xi}(t+7r)|_2^2\leq E_{11},\eqno{(5.98)}$$ where
\begin{eqnarray*}
 E_{11}&=&c[\frac{E_1}{r}+E_1+E_2+E_4+E_4^2+E_2^{\frac{1}{2}}(E_6^{\frac{1}{2}}+E_8)+(E_4^2+E_2^2)E_9\\
&+&E_{10}+|Q_1|_2^2+|Q_{1\xi}|_2^2+|Q_2|_2^2+|Q_{2\xi}|_2^2+|Q_1|_{\xi=1}|_2^2+|Q_2|_{\xi=1}|_2^2]\\
&\cdot&\exp c(1+E_2^2+E_4^2+E_6^2+E_8^4 ).
\end{eqnarray*}
From (5.97) and (5.98), we have
$$c_1\int_{t+7r}^{t+8r}(\|T_{\xi}\|^2+\|q_{\xi}\|^2)\leq E_{11}^2+2E_{11}+E_1=E_{12}.\eqno{(5.99)}$$
By Gronwall inequality, from (5.97) we obtain
$$|T_{\xi}(t)|_2^2+|q_{\xi}(t)|_2^2\leq
C_5,\eqno{(5.100)}$$ where $C_5=C_5(\|U_0\|,\ \|Q_1\|_1,\
\|Q_2\|_1)>0$ and $0\leq t<7r $.

\noindent{\bf $H^1$ estimates about $v,\ T,\ q$} \quad Taking the
inner product of equation (2.11) with $-\triangle v$ in
$L^2(\Omega)\times L^2(\Omega)$, we get
\begin{eqnarray*}
& &\frac{1}{2}\frac{d
\displaystyle\int_{\Omega}(|\nabla_{e_{\theta}}v|^2+|\nabla_{e_{\varphi}}v|^2+|v|^2)}{d
t}+\frac{1}{Re_1}\displaystyle\int_{\Omega}|\triangle v|_2^2
\\&+&\frac{1}{Re_2}\displaystyle\int_{\Omega}(|\nabla_{e_{\theta}}v_{\xi}|^2
+|\nabla_{e_{\varphi}}v_{\xi}|^2+|v_{\xi}|^2)\\
&=&\displaystyle\int_{\Omega}(\nabla_{v}v+W(v)v_{\xi})\cdot\triangle
v+\displaystyle\int_{\Omega}[\int_{\xi}^1\frac{b
P}{p}\hbox{grad}((1+aq)T)d
\xi^{'}]\cdot\triangle v\\
&+&\displaystyle\int_{\Omega}(\frac{f}{R_0}k\times v+\hbox{grad}
\Phi_s)\cdot\triangle v.
\qquad\qquad\qquad\qquad\quad\quad\quad\quad\quad\quad\quad
(5.101)
\end{eqnarray*}
By H\"older inequality, (4.16) and Young inequality, we have
\begin{eqnarray*}
&&|\displaystyle\int_{\Omega}\nabla_{v}v\cdot\triangle v|\\
&\leq&\displaystyle\int_{\Omega}|v|(|\nabla_{e_{\theta}}v|^2+|\nabla_{e_{\varphi}}v|^2)^{\frac{1}{2}}|\triangle
v|\\
&\leq&c|v|_4^2(\displaystyle\int_{\Omega}(|\nabla_{e_{\theta}}v|^2+|\nabla_{e_{\varphi}}v|^2)^2)^{\frac{1}{2}}
+\varepsilon|\triangle v|^2_2\\
&\leq&c|v|_4^2(\displaystyle\int_{\Omega}(|\nabla_{e_{\theta}}v|^2
+|\nabla_{e_{\varphi}}v|^2))^{\frac{1}{4}}(\displaystyle\int_{\Omega}(|\nabla_{e_{\theta}}v|^2
+|\nabla_{e_{\varphi}}v|^2\\
&+&|\nabla_{e_{\theta}}v_{\xi}|^2
+|\nabla_{e_{\varphi}}v_{\xi}|^2+|\triangle
v|^2))^{\frac{3}{4}}+\varepsilon|\triangle v|^2_2\\
&\leq&c(|v|_4^8+|v|_4^2)\displaystyle\int_{\Omega}(|\nabla_{e_{\theta}}v|^2+|\nabla_{e_{\varphi}}v|^2)
+2\varepsilon(|\triangle
v|^2_2+\displaystyle\int_{\Omega}(|\nabla_{e_{\theta}}v_\xi|^2+|\nabla_{e_{\varphi}}v_\xi|^2)).
\ (5.102)
\end{eqnarray*}
By H\"older inequality, Minkowski inequality, (4.13) and Young
inequality, we obtain
\begin{eqnarray*}
&&|\displaystyle\int_{\Omega}W(v)v_{\xi}\cdot\triangle
v|\\
&\leq&\displaystyle\int_{S^2}[\int_{0}^{1}(|\nabla_{e_{\theta}}v|^2
+|\nabla_{e_{\varphi}}v|^2)^{\frac{1}{2}}d\xi\int_{0}^{1}|v_{\xi}||\triangle
v|d\xi]\\
&\leq&\displaystyle\int_{S^2}[(\int_{0}^{1}(|\nabla_{e_{\theta}}v|^2
+|\nabla_{e_{\varphi}}v|^2)d\xi)^{\frac{1}{2}}(\int_{0}^{1}|v_{\xi}|^2d\xi)^{\frac{1}{2}}(\int_{0}^{1}|\triangle
v|^2d\xi)^{\frac{1}{2}}]\\
&\leq&[\displaystyle\int_{S^2}(\int_{0}^{1}(|\nabla_{e_{\theta}}v|^2
+|\nabla_{e_{\varphi}}v|^2)d\xi)^2]^{\frac{1}{4}}[\displaystyle\int_{S^2}
(\int_{0}^{1}|v_{\xi}|^2d\xi)^2]^{\frac{1}{4}}|\triangle
v|_2\\
&\leq&c[\displaystyle\int_{0}^{1}(\int_{S^2}(|\nabla_{e_{\theta}}v|^2
+|\nabla_{e_{\varphi}}v|^2)^2)^{\frac{1}{2}}d\xi][\displaystyle\int_{0}^{1}
(\int_{S^2}|v_{\xi}|^4)^{\frac{1}{2}}d\xi]+\varepsilon|\triangle|_2^2\\
&\leq&c[\displaystyle\int_{0}^{1}(\int_{S^2}(|\nabla_{e_{\theta}}v|^2
+|\nabla_{e_{\varphi}}v|^2))^{\frac{1}{2}}\cdot(\int_{S^2}(|\nabla_{e_{\theta}}v|^2
+|\nabla_{e_{\varphi}}v|^2+|\triangle
v|^2))^{\frac{1}{2}}d\xi]\\
&\cdot&c[\displaystyle\int_{0}^{1}(\int_{S^2}|v_{\xi}|^2)^{\frac{1}{2}}\cdot(\int_{S^2}
(|\nabla_{e_{\theta}}v_{\xi}|^2+|\nabla_{e_{\varphi}}v_{\xi}|^2+|
v_{\xi}|^2))^{\frac{1}{2}}d\xi]+\varepsilon|\triangle
v|_2^2\\
&\leq&c[\displaystyle\int_{\Omega}(|\nabla_{e_{\theta}}v|^2+|\nabla_{e_{\varphi}}v|^2)
+(\int_{\Omega}(|\nabla_{e_{\theta}}v|^2+|\nabla_{e_{\varphi}}v|^2))^{\frac{1}{2}}|\triangle
v|_2]\\
&\cdot&
|v_{\xi}|_2[|v_{\xi}|_2+(\int_{\Omega}(|\nabla_{e_{\theta}}v_\xi|^2
+|\nabla_{e_{\varphi}}v_\xi|^2))^{\frac{1}{2}}]+\varepsilon|\triangle
v|_2^2\\
&\leq&2\varepsilon|\triangle
v|_2^2+c[2|v_{\xi}|^2_2+|v_{\xi}|^4_2+(|v_{\xi}|^2_2+1)\displaystyle\int_{\Omega}(|\nabla_{e_{\theta}}v_{\xi}|^2
+|\nabla_{e_{\varphi}}v_{\xi}|^2)]\\
&\cdot&\int_{\Omega}(|\nabla_{e_{\theta}}v|^2
+|\nabla_{e_{\varphi}}v|^2).
\quad\quad\quad\quad\quad\quad\quad\quad\quad\quad\quad\quad\quad\quad\quad\quad\quad\quad(5.103)
\end{eqnarray*}
By H\"older inequality, Young inequality, Minkowski inequality,
and (4.13), we have
\begin{eqnarray*}
& &|\displaystyle\int_{\Omega}\int_{\xi}^{1}\frac{b
P}{p}\hbox{grad}[(1+a
q)T]d\xi^{'}\cdot\triangle v|\\
&\leq&c\displaystyle\int_{\Omega}(\int_{0}^{1}|\nabla
T|d\xi|\triangle v|+\int^1_0|\nabla  T|qd\xi|\triangle
v|+\int_{0}^{1}|T||\nabla q|d\xi|\triangle v|)\\
&\leq&c|\nabla T|_2|\triangle
v|_2+c\displaystyle\int_{\Omega}[(\int_{0}^{1}|q|^2d\xi)^{\frac{1}{2}}(\int_{0}^{1}|\nabla
T|^2d\xi)^{\frac{1}{2}}|\triangle
v|]\\
&+&c\displaystyle\int_{\Omega}[(\int_{0}^{1}|T|^2d\xi)^{\frac{1}{2}}(\int_{0}^{1}|\nabla
q|^2d\xi)^{\frac{1}{2}}|\triangle v|]\\
&\leq&c|\nabla T|_2|\triangle
v|_2+[(\displaystyle\int_{\Omega}(\int_{0}^{1}|q|^2d\xi)^2)^{\frac{1}{4}}(\displaystyle\int_{\Omega}(\int_{0}^{1}|\nabla
T|^2d\xi)^2)^{\frac{1}{4}}\\
&+&(\displaystyle\int_{\Omega}(\int_{0}^{1}|T|^2d\xi)^2)^{\frac{1}{4}}(\displaystyle\int_{\Omega}(\int_{0}^{1}|\nabla
q|^2d\xi)^2)^{\frac{1}{4}}]|\triangle
v|_2\\
&\leq&c(\displaystyle\int_{\Omega}(\int_{0}^{1}|q|^2d\xi)^2)^{\frac{1}{2}}(\displaystyle\int_{\Omega}(\int_{0}^{1}|\nabla
T|^2d\xi)^2)^{\frac{1}{2}}\\&+&c(\displaystyle\int_{\Omega}(\int_{0}^{1}|T|^2d\xi)^2)^{\frac{1}{2}}
(\displaystyle\int_{\Omega}(\int_{0}^{1}|\nabla
q|^2d\xi)^2)^{\frac{1}{2}}+c|\nabla T|_2^2+\varepsilon|\triangle
v|_2^2\\
&\leq&c(\displaystyle\int_{0}^{1}(\int_{S^2}|q|^4)^{\frac{1}{2}}d\xi)(\displaystyle\int_{0}^{1}(\int_{S^2}|\nabla
T|^4)^{\frac{1}{2}}d\xi)\\&+&c(\displaystyle\int_{0}^{1}(\int_{S^2}|T|^4)^{\frac{1}{2}}d\xi)
(\displaystyle\int_{0}^{1}(\int_{S^2}|\nabla
q|^4)^{\frac{1}{2}}d\xi)+c|\nabla T|_2^2+\varepsilon|\triangle
v|_2^2\\
&\leq&c|q|_4^2\displaystyle\int_{0}^{1}[\|\nabla
T\|_{L^2(S^2)}(\|\nabla T\|_{L^2(S^2)}^2+\|\triangle
T\|_{L^2(S^2)}^2)^{\frac{1}{2}}]d\xi\\
&+&c|T|_4^2\displaystyle\int_{0}^{1}[\|\nabla
q\|_{L^2(S^2)}(\|\nabla q\|_{L^2(S^2)}^2+\|\triangle
q\|_{L^2(S^2)}^2)^{\frac{1}{2}}]d\xi+c|\nabla
T|_2^2+\varepsilon|\triangle
v|_2^2\\
&\leq&c|q|_4^2|\nabla T|_2(|\nabla T|_2^2+|\triangle
T|_2^2)^{\frac{1}{2}}+c|T|_4^2|\nabla q|_2(|\nabla
q|_2^2+|\triangle q|_2^2)^{\frac{1}{2}}\\&+&c|\nabla
T|_2^2+\varepsilon|\triangle
v|_2^2\\
&\leq&c|q|_4^2|\nabla T|_2^2+c|T|_4^2|\nabla
q|_2^2+c|q|_4^4|\nabla T|_2^2+c|T|_4^4|\nabla q|_2^2+c|\nabla
T|_2^2\\&+&\varepsilon|\triangle v|_2^2+\varepsilon|\triangle
T|_2^2+\varepsilon|\triangle q|_2^2.\qquad \qquad\quad\quad
\qquad\qquad\qquad\qquad\quad\quad (5.104)
\end{eqnarray*}
By using Lemma 4.1 and (2.19), we get
$$\displaystyle\int_{\Omega}\hbox{grad} \Phi_s\cdot\triangle v=0.\eqno(5.105)$$
$(\frac{f}{R_0}k\times v)\cdot\triangle v=0$ implies
$$\displaystyle\int_{\Omega}(\frac{f}{R_0}k\times v)\cdot\triangle v=0.\eqno(5.106)$$
We derive from (5.101)-(5.106)
\begin{eqnarray*}
& &\frac{1}{2}\frac{d
\displaystyle\int_{\Omega}(|\nabla_{e_{\theta}}v|^2+|\nabla_{e_{\varphi}}v|^2+|v|^2)}{d
t}+\frac{1}{Re_1}|\triangle v|_2^2
\\&+&\frac{1}{Re_2}\displaystyle\int_{\Omega}(|\nabla_{e_{\theta}}v_{\xi}|^2
+|\nabla_{e_{\varphi}}v_{\xi}|^2+|v_{\xi}|^2)\\
&\leq&c[|v|_4^8+|v|_4^2+2|v_{\xi}|^2_2+|v_{\xi}|^4_2+(|v_{\xi}|^2_2+1)
\|v_{\xi}\|^2
]\displaystyle\int_{\Omega}(|\nabla_{e_{\theta}}v|^2+|\nabla_{e_{\varphi}}v|^2)
\\&+&c(1+|q|_4^2+|q|_4^4)|\nabla T|_2^2+c(|T|_4^2+|T|_4^4)|\nabla
q|_2^2+2\varepsilon\displaystyle\int_{\Omega}(|\nabla_{e_{\theta}}v_{\xi}|^2
+|\nabla_{e_{\varphi}}v_{\xi}|^2)\\
&+&5\varepsilon|\triangle v|_2^2+\varepsilon|\triangle
T|_2^2+\varepsilon|\triangle q|_2^2+c|\nabla T|_2^2.\quad\qquad
\qquad\qquad\qquad\quad\quad(5.107)
\end{eqnarray*}

\bigskip
We take the inner product of equation (2.12) with $-\triangle T$
in $L^2(\Omega)$ and obtain
\begin{eqnarray*}
& &\frac{1}{2}\frac{d |\nabla T|_2^2}{d
t}+\frac{1}{Rt_1}|\triangle T|_2^2+\frac{1}{Rt_2}(|\nabla
T_{\xi}|_2^2+\alpha_s |\nabla T|_{\xi=1}|_2^2)\\
&=&\displaystyle\int_{\Omega}(\nabla_{v}T+W(v)\frac{\partial
T}{\partial \xi})\triangle
T-\displaystyle\int_{\Omega}\frac{bP}{p}(1+aq)W(v)\triangle T
-\displaystyle\int_{\Omega}Q_1\triangle T. (5.108)
\end{eqnarray*}
Similarly to (5.102), we get
$$|\displaystyle\int_{\Omega}\triangle
T\nabla_{v}T|\leq c(|v|_4^8+|v|_4^2)|\nabla
T|_2^2+2\varepsilon(|\triangle T|_2^2+|\nabla
T_{\xi}|_2^2).\eqno(5.109)$$ We derive as (5.105)
\begin{eqnarray*}
|\displaystyle\int_{\Omega}W(v)T_{\xi}\triangle
T|&\leq&\varepsilon|\triangle T|_2^2+\varepsilon|\triangle v|_2^2
+c[2|T_{\xi}|^2_2+|T_{\xi}|^4_2\\
&+&(|T_{\xi}|^2_2+1)\displaystyle\int_{\Omega}|\nabla T_{\xi}|^2
]\displaystyle\int_{\Omega}(|\nabla_{e_{\theta}}v|^2+|\nabla_{e_{\varphi}}v|^2)
. \ \quad(5.110)
\end{eqnarray*}
By H\"older inequality, Minkowsky inequality, (4.13) and Young
inequality, we obtain
\begin{eqnarray*}
& &|\displaystyle\int_{\Omega}\frac{b P}{p}(1+a q)W(v)\triangle
T|\\
&\leq&c\displaystyle\int_{\Omega}(\int_{0}^1|\hbox{div}v|d\xi)|\triangle
T|+c\displaystyle\int_{\Omega}|q|(\int_{0}^1|\hbox{div}v|d\xi)|\triangle
T|\\
&\leq&c(\displaystyle\int_{\Omega}|\hbox{div}v|^2)^{\frac{1}{2}}|\triangle
T|_2+c|q|_4(\displaystyle\int_{\Omega}(\int_{0}^{1}|\hbox{div}v|^2d\xi)^2)^{\frac{1}{4}}|\triangle
T|_2\\
&\leq&\varepsilon|\triangle
T|_2^2+c\displaystyle\int_{\Omega}(|\nabla_{e_{\theta}}v|^2+|\nabla_{e_{\varphi}}v|^2)\\
&+&c|q|_4^2\displaystyle\int_{0}^{1}[(\int_{S^2}(|\nabla_{e_{\theta}}v|^2
+|\nabla_{e_{\varphi}}v|^2))^{\frac{1}{2}}(\int_{S^2}(|\nabla_{e_{\theta}}v|^2+|\nabla_{e_{\varphi}}v|^2
+|\triangle v|^2))^{\frac{1}{2}}]d\xi\\
&\leq&\varepsilon|\triangle T|_2^2+\varepsilon|\triangle
v|_2^2+c(1+|q|_4^2+|q|_4^4)\displaystyle\int_{\Omega}(|\nabla_{e_{\theta}}v|^2
+|\nabla_{e_{\varphi}}v|^2).\quad\quad(5.111)
\end{eqnarray*}
By H\"older inequality and Young inequality, we have
$$|\displaystyle\int_{\Omega}Q_1\triangle T|\leq c|Q_1|_2^2+\varepsilon|\triangle
T|_2^2.\eqno(5.112)$$ From (5.108)-(5.112), we obtain
\begin{eqnarray*}
&&\frac{1}{2}\frac{d |\nabla T|_2^2}{d t}+\frac{1}{Rt_1}|\triangle
T|_2^2+\frac{1}{Rt_2}(|\nabla T_{\xi}|_2^2+\alpha_s |\nabla
T|_{\xi=1}|_2^2)
\\
&\leq& 5\varepsilon|\triangle T|_2^2+2\varepsilon|\nabla
T_{\xi}|_2^2+2\varepsilon|\triangle v|_2^2+
c[1+|q|_4^2+|q|_4^4+2|T_{\xi}|^2_2\\
&+&|T_{\xi}|^4_2+(|T_{\xi}|^2_2+1)\displaystyle\int_{\Omega}|\nabla
T_{\xi}|^2]
\displaystyle\int_{\Omega}(|\nabla_{e_{\theta}}v|^2+|\nabla_{e_{\varphi}}v|^2)\\
&+&c(|v|_4^2+|v|_4^8)|\nabla
T|_2^2+c|Q_1|_2^2.\quad\quad\quad\quad\quad\quad
\quad\quad\quad\quad\quad\quad\quad\quad\quad(5.113)
\end{eqnarray*}
By taking the inner product of equation (2.13) with $-\triangle q$
in $L^2(\Omega)$, we have
\begin{eqnarray*}
& &\frac{1}{2}\frac{d |\nabla q|_2^2}{d
t}+\frac{1}{Rq_1}|\triangle q|_2^2+\frac{1}{Rq_2}(|\nabla
q_{\xi}|_2^2+\beta_s |\nabla
q|_{\xi=1}|_2^2)\\
&=&\displaystyle\int_{\Omega}(\nabla_{v}q+W(v)\frac{\partial
q}{\partial \xi})\triangle
q-\displaystyle\int_{\Omega}Q_2\triangle q. \quad\quad \
\quad\quad\quad\quad\quad\quad\quad\quad(5.114)
\end{eqnarray*}
Similarly to (5.113), we derive from (5.114)
\begin{eqnarray*}
& &\frac{1}{2}\frac{d |\nabla q|_2^2}{d
t}+\frac{1}{Rq_1}|\triangle q|_2^2+\frac{1}{Rq_2}(|\nabla
q_{\xi}|_2^2+\beta_s |\nabla
q|_{\xi=1}|_2^2)\\
&\leq&4\varepsilon|\triangle q|_2^2+2\varepsilon|\nabla
q_{\xi}|_2^2+\varepsilon|\triangle
v|_2^2+c(|v|_4^2+|v|_4^8)|\nabla
q|_2^2+c[2|q_{\xi}|^2_2+|q_{\xi}|^4_2\\
&+&(|q_{\xi}|^2_2+1)\displaystyle\int_{\Omega}|\nabla
q_{\xi}|^2]\displaystyle\int_{\Omega}(|\nabla_{e_{\theta}}v|^2
+|\nabla_{e_{\varphi}}v|^2)+c|Q_2|_2^2.\quad\quad\quad\quad(5.115)
\end{eqnarray*}

From (5.107), (5.113) and (5.115), choosing $\varepsilon$ small
enough, we obtain
\begin{eqnarray*}
& &\frac{d
[\displaystyle\int_{\Omega}(|\nabla_{e_{\theta}}v|^2+|\nabla_{e_{\varphi}}v|^2+|v|^2)+|\nabla
T|_2^2+|\nabla q|_2^2]}{d
t}\\
&+&\frac{1}{Re_1}|\triangle v|_2^2+\frac{1}{Rt_1}|\triangle
T|_2^2+\frac{1}{Rq_1}|\triangle
q|_2^2+\frac{1}{Re_2}\displaystyle\int_{\Omega}(|\nabla_{e_{\theta}}v_{\xi}|^2
+|\nabla_{e_{\varphi}}v_{\xi}|^2+|v_{\xi}|^2)\\
&+&\frac{1}{Rt_2}(|\nabla T_{\xi}|_2^2+\alpha_s |\nabla
T|_{\xi=1}|_2^2)+\frac{1}{Rq_2}(|\nabla q_{\xi}|_2^2+\beta_s
|\nabla
q|_{\xi=1}|_2^2)\\
&\leq&c[1+|q|_4^2+|q|_4^4+|T|_4^2+|T|_4^4+|v|_4^2+|v|_4^8+2|v_{\xi}|^2_2+|v_{\xi}|^4_2+(|v_{\xi}|^2_2+1)
\|v_{\xi}\|^2\\
&+&2|T_{\xi}|^2_2
+|T_{\xi}|^4_2+(|T_{\xi}|^2_2+1)\displaystyle\int_{\Omega}|\nabla
T_{\xi}|^2+2|q_{\xi}|^2_2+|q_{\xi}|^4_2+(|q_{\xi}|^2_2+1)\displaystyle\int_{\Omega}|\nabla
q_{\xi}|^2]\\
&\cdot&[\displaystyle\int_{\Omega}(|\nabla_{e_{\theta}}v|^2
+|\nabla_{e_{\varphi}}v|^2)+|\nabla T|_2^2+|\nabla
q|_2^2]+c|Q_1|_2^2+c|Q_2|_2^2. \quad\quad\ (5.116)
\end{eqnarray*}
By Lemma 4.6, (5.18), (5.25), (5.39), (5.64), (5.73), (5.82),
(5.83), (5.98) and (5.99) , we get
$$|\nabla_{e_{\theta}}v(t+8r)|_2^2+|\nabla_{e_{\varphi}}v(t+8r)|_2^2+|\nabla T(t+8r)|_2^2
+|\nabla q(t+8r)|_2^2\leq E_{13},\eqno{(5.117)}$$ where
\begin{eqnarray*}
 E_{13}&=&c(\frac{E_1}{r}+|Q_1|_2^2
+|Q_2|_2^2) \cdot\exp
c[1+E_2+E_4+E_6^2+E_8^4+E_9^2\\
&+&(E_9+1)E_{10}+E_{11}^2+(E_{11}+1)E_{12}].
\end{eqnarray*}
By Gronwall inequality, from (5.116) we obtain
$$|\nabla_{e_{\theta}}v(t)|_2^2+|\nabla_{e_{\varphi}}v(t)|_2^2+|\nabla T(t)|_2^2
+|\nabla q(t)|_2^2\leq C_6,\eqno{(5.118)}$$ where
$C_6=C_6(\|U_0\|,\ \|Q_1\|,\ \|Q_2\|)>0$ and $0\leq t<8r $.

\subsection{The global existence of strong solutions}

\smallskip
\noindent{\bf Proof of Theorem 3.1.} \quad By Proposition 5.3, we
can use the method of contradiction to prove Theorem 3.1. Indeed,
let $U$ be a strong solution to the system (2.11)-(2.17) on the
maximal interval $[0,\mathcal{T}_*]$. If $\mathcal{T}_*<+\infty$,
then
$$\displaystyle\limsup_{t\rightarrow
\mathcal{T}_*^-}\|U\|=+\infty,$$  which is impossible from (5.18),
(5.82), (5.84), (5.98), (5.100), (5.117), (5.118). The proof is
complete.

\smallskip
\section{The uniqueness of strong solutions}

\smallskip
\noindent{\bf Proof of Theorem 3.2} \quad Let $(v_1, T_1, q_1)$
and $(v_2, T_2, q_2)$ be two strong solutions of (2.11)-(2.17) on
the time interval $[0,\mathcal{T}]$ with corresponding
geopotentials $\Phi_{s_1}$, $\Phi_{s_2}$, and initial data
$((v_0)_1, (T_0)_1, (q_0)_1)$, $((v_0)_2, (T_0)_2, (q_0)_2)$,
respectively.

Define $v=v_1-v_2$, $T=T_1-T_2$, $q=q_1-q_2$,
$\Phi_s=\Phi_{s_1}-\Phi_{s_2}$. Then $v$, $T$, $q$, $\Phi_s$
satisfy the following system
$$\frac{\partial v}{\partial t}-\frac{1}{Re_1}\triangle v-\frac{1}{Re_2}\frac{\partial^2 v}{\partial \xi^2}
+\nabla_{v_1}v+\nabla_{v}v_2+W(v_1)\frac{\partial v}{\partial
\xi}+ W(v)\frac{\partial v_2}{\partial \xi}+\frac{f}{R_0}k\times
v$$
$$+\hbox{grad} \Phi_s+\displaystyle\int_{\xi}^1\frac{bP}{p}\hbox{grad}Td\xi^{'}
+\displaystyle\int_{\xi}^1\frac{abP}{p}\hbox{grad} (q_1
T)d\xi^{'}+ \displaystyle\int_{\xi}^1\frac{abP}{p}\hbox{grad} (q
T_2)d\xi^{'}=0, \eqno(6.1)$$
$$\frac{\partial T}{\partial t}-\frac{1}{Rt_1}\triangle T-\frac{1}{Rt_2}\frac{\partial^2 T}{\partial \xi^2}
+\nabla_{v_1}T+\nabla_{v}T_2+W(v_1)\frac{\partial T}{\partial \xi}+
W(v)\frac{\partial T_2}{\partial \xi}-\frac{bP}{p}W(v)$$
$$-\frac{abP}{p}q_1W(v)-\frac{abP}{p}qW(v_2)=0, \eqno(6.2)$$
$$\frac{\partial q}{\partial t}-\frac{1}{Rq_1}\triangle q-\frac{1}{Rq_2}\frac{\partial^2 q}{\partial \xi^2}
+\nabla_{v_1}q+\nabla_{v}q_2+W(v_1)\frac{\partial q}{\partial
\xi}+ W(v)\frac{\partial q_2}{\partial \xi}=0,
\quad\quad\quad\quad\eqno(6.3)$$
$$v|_{t=0}=(v_0)_1-(v_0)_2,\quad\quad\quad\quad\quad\quad\quad\quad\quad\quad\quad\quad
\quad\quad\quad\quad\quad\quad\quad\quad\quad\quad\quad\eqno(6.4)$$
$$T|_{t=0}=(T_0)_1-(T_0)_2,\quad\quad\quad\quad\quad\quad\quad\quad\quad\quad\quad\quad
\quad\quad\quad\quad\quad\quad\quad\quad\quad\quad\quad\eqno(6.5)$$
$$q|_{t=0}=(q_0)_1-(q_0)_2,\quad\quad\quad\quad\quad\quad\quad\quad\quad\quad\quad\quad
\quad\quad\quad\quad\quad\quad\quad\quad\quad\quad\quad\eqno(6.6)$$
$$\xi=1:\ \frac{\partial v}{\partial \xi}=0,\ \frac{\partial T}{\partial \xi}=
-\alpha_s T,\ \frac{\partial q}{\partial \xi}= -\beta_s q,
\quad\quad\quad\quad\quad\quad\quad\quad\quad\quad\quad\quad\quad
$$
$$\xi=0:\ \frac{\partial v}{\partial \xi}=0,\ \frac{\partial T}{\partial \xi}=0,\
\frac{\partial q}{\partial \xi}=
0.\quad\quad\quad\quad\quad\quad\quad\quad\quad\quad\quad\quad\quad\quad\quad$$
We take the inner product of equation (6.1) with $v$ in
$L^2(\Omega)\times L^2(\Omega)$ and obtain
\begin{eqnarray*}
& &\frac{1}{2}\frac{d |v|_2^2}{d
t}+\frac{1}{Re_1}\displaystyle\int_{\Omega}(|\nabla_{e_{\theta}}v|^2+|\nabla_{e_{\varphi}}v|^2+|v|^2)
+\frac{1}{Re_2}\displaystyle\int_{\Omega}|v_{\xi}|^2\\
&=&-\displaystyle\int_{\Omega}(\nabla_{v_1}v+W(v_1)\frac{\partial
v}{\partial \xi})\cdot
v-\displaystyle\int_{\Omega}v\cdot\nabla_{v}v_2
-\displaystyle\int_{\Omega}W(v)\frac{\partial v_2}{\partial
\xi}\cdot v\\&-&\displaystyle\int_{\Omega}(\frac{f}{R_0}k\times
v+\hbox{grad} \Phi_s)\cdot
v-\displaystyle\int_{\Omega}(\displaystyle\int_{\xi}^1\frac{bP}{p}\hbox{grad}
Td\xi^{'})\cdot
v\\&-&\displaystyle\int_{\Omega}(\displaystyle\int_{\xi}^1\frac{abP}{p}\hbox{grad}
(q_1 T)d\xi^{'})\cdot
v-\displaystyle\int_{\Omega}(\displaystyle\int_{\xi}^1\frac{abP}{p}\hbox{grad}
(q T_2)d\xi^{'})\cdot v.\quad\qquad(6.7)
\end{eqnarray*}
By Lemma 4.4, we have
$$\displaystyle\int_{\Omega}(\nabla_{v_1}v+W(v_1)\frac{\partial
v}{\partial \xi})\cdot v=0. \eqno(6.8)$$ Using Lemma 4.3, H\"older
inequality, Young inequality and (4.16), we get
\begin{eqnarray*}
|\displaystyle\int_{\Omega}v\cdot\nabla_{v}v_2
|&=&|\displaystyle\int_{\Omega}(v_2\cdot\nabla_{v}v+v_2\cdot v
\hbox{div}
v)|\\
&\leq&c\displaystyle\int_{\Omega}|v||v_2|(|\nabla_{e_{\theta}}v|^2+|\nabla_{e_{\varphi}}v|^2)^{\frac{1}{2}}\\
&\leq&\varepsilon\displaystyle\int_{\Omega}(|\nabla_{e_{\theta}}v|^2+|\nabla_{e_{\varphi}}v|^2)+c|v|_4^2|v_2|_4^2\\
&\leq&\varepsilon\displaystyle\int_{\Omega}(|\nabla_{e_{\theta}}v|^2
+|\nabla_{e_{\varphi}}v|^2)+c|v|_2^{\frac{1}{2}}|v_2|_4^2\|v\|^{\frac{3}{2}}\\
&\leq&2\varepsilon\|v\|^2+c|v_2|_4^8|v|_2^2.\qquad\qquad\qquad\qquad\qquad\qquad\
\ \ \ (6.9)
\end{eqnarray*}
By H\"older inequality, Young inequality, Minkowski inequality and
(4.13), we obtain
\begin{eqnarray*}
& &|\displaystyle\int_{\Omega}W(v)\frac{\partial v_2}{\partial
\xi}\cdot v|\\
&\leq& \displaystyle\int_{S^2}[\int_{0}^1(|\nabla_{e_{\theta}}v|^2
+|\nabla_{e_{\varphi}}v|^2)^{\frac{1}{2}}d\xi\int_{0}^1|v_{2\xi}||v|d\xi]
\\&\leq&(\displaystyle\int_{\Omega}(|\nabla_{e_{\theta}}v|^2
+|\nabla_{e_{\varphi}}v|^2))^{\frac{1}{2}}(\displaystyle\int_{S^2}
(\int_{0}^1|v_{2\xi}|^2d\xi)^2)^{\frac{1}{4}}(\displaystyle\int_{S^2}(\int_{0}^1|v|^2d\xi)^2)^{\frac{1}{4}}\\
&\leq&\varepsilon\displaystyle\int_{\Omega}(|\nabla_{e_{\theta}}v|^2+|\nabla_{e_{\varphi}}v|^2)
+c\displaystyle\int_{0}^1(\int_{S^2}|v_{2\xi}|^4)^{\frac{1}{2}}d\xi\displaystyle\int_{0}^1
(\int_{S^2}|v|^4)^{\frac{1}{2}}d\xi\\
&\leq&\varepsilon\displaystyle\int_{\Omega}(|\nabla_{e_{\theta}}v|^2+|\nabla_{e_{\varphi}}v|^2)
+c\displaystyle\int_{0}^1[\|v_{2\xi}\|_{L^2(S^2)}(\int_{S^2}(|\nabla_{e_{\theta}}v_{2\xi}|^2
+|\nabla_{e_{\varphi}}v_{2\xi}|^2))^{\frac{1}{2}}]d\xi\\
&\cdot&\displaystyle\int_{0}^1
[\|v\|_{L^2(S^2)}(\int_{S^2}(|\nabla_{e_{\theta}}v|^2+|\nabla_{e_{\varphi}}v|^2+|v|^2))^{\frac{1}{2}}]d\xi\\
&\leq&2\varepsilon\|v\|^2+c[(|v_{2\xi}|_2^2+1)\displaystyle\int_{\Omega}(|\nabla_{e_{\theta}}v_{2\xi}|^2
+|\nabla_{e_{\varphi}}v_{2\xi}|^2)+|v_{2\xi}|_2^2
]|v|_2^2.\quad\quad\quad \ (6.10)
\end{eqnarray*}
By Lemma 4.1, H\"older inequality, Young inequality, Minkowski
inequality and (4.13), we have
\begin{eqnarray*}
&
&|\displaystyle\int_{\Omega}(\displaystyle\int_{\xi}^1\frac{abP}{p}\hbox{grad}
(q T_2)d\xi^{'})\cdot v|\\
&=&|\displaystyle\int_{\Omega}(\displaystyle\int_{\xi}^1\frac{abP}{p}qT_2d\xi^{'})\cdot
\hbox{div} v|\\
&\leq&c\displaystyle\int_{S^2}[\int_{0}^1|q||T_2|d\xi\int_{0}^1(|\nabla_{e_{\theta}}v|^2
+|\nabla_{e_{\varphi}}v|^2)^{\frac{1}{2}}d\xi]\\
&\leq&c(\displaystyle\int_{S^2}(\int_{0}^1|q|^2d\xi)^2)^{\frac{1}{2}}
(\displaystyle\int_{S^2}(\int_{0}^1|T_2|^2d\xi)^2)^{\frac{1}{2}}
+\varepsilon\displaystyle\int_{\Omega}(|\nabla_{e_{\theta}}v|^2+|\nabla_{e_{\varphi}}v|^2)\\
&\leq&c\displaystyle\int_{0}^1(\int_{S^2}|q|^4)^{\frac{1}{2}}d\xi\displaystyle\int_{0}^1
(\int_{S^2}|T_2|^4)^{\frac{1}{2}}d\xi+\varepsilon\displaystyle\int_{\Omega}(|\nabla_{e_{\theta}}v|^2
+|\nabla_{e_{\varphi}}v|^2)\\
&\leq&c|q|_2(|\nabla q|_2+|q|_2)|T_2|_4^2+\varepsilon\displaystyle\int_{\Omega}(|\nabla_{e_{\theta}}v|^2
+|\nabla_{e_{\varphi}}v|^2)\\
&\leq&c(|T_2|_4^2+|T_2|_4^4)|q|_2^2+\varepsilon|\nabla
q|_2^2+\varepsilon\displaystyle\int_{\Omega}(|\nabla_{e_{\theta}}v|^2+|\nabla_{e_{\varphi}}v|^2).\qquad
\quad\quad (6.11)
\end{eqnarray*}
By $(\frac{f}{R_0}k\times v)\cdot v=0$ and Lemma 4.1, we get
$$\displaystyle\int_{\Omega}(\frac{f}{R_0}k\times v+\hbox{grad} \Phi_s)\cdot v=0.\eqno(6.12)$$
We derive from (6.7)-(6.12)
\begin{eqnarray*}
& &\frac{1}{2}\frac{d |v|_2^2}{d
t}+\frac{1}{Re_1}\displaystyle\int_{\Omega}(|\nabla_{e_{\theta}}v|^2+|\nabla_{e_{\varphi}}v|^2+|v|^2)
+\frac{1}{Re_2}\displaystyle\int_{\Omega}|v_{\xi}|^2\\
&\leq&5\varepsilon\|v\|^2+\varepsilon|\nabla
q|_2^2+c(|T_2|_4^2+|T_2|_4^4)|q|_2^2\\
&+&c[|v_2|_4^8+|v_{2\xi}|_2^2+(|v_{2\xi}|_2^2+1)\displaystyle\int_{\Omega}(|\nabla_{e_{\theta}}v_{2\xi}|^2
+|\nabla_{e_{\varphi}}v_{2\xi}|^2)]|v|_2^2\\
&-&\displaystyle\int_{\Omega}(\displaystyle\int_{\xi}^1\frac{b
P}{p}\hbox{grad} Td\xi^{'})\cdot
v-\displaystyle\int_{\Omega}(\displaystyle\int_{\xi}^1\frac{ab
P}{p}\hbox{grad}(q_1T)d\xi^{'})\cdot v.\quad\qquad\qquad(6.13)
\end{eqnarray*}
By taking the inner product of equation (6.2) with $T$ in
$L^2(\Omega)$ and equation (6.3) with $q$ in $L^2(\Omega)$, we
obtain
\begin{eqnarray*}
& &\frac{1}{2}\frac{d |T|_2^2}{d
t}+\frac{1}{Rt_1}\displaystyle\int_{\Omega}|\nabla T|^2
+\frac{1}{Rt_2}\displaystyle\int_{\Omega}|T_{\xi}|^2+\frac{\alpha_s}{Rt_2}|T|_{\xi=1}|_2^2\\
&=&-\displaystyle\int_{\Omega}(\nabla_{v_1}T+W(v_1)\frac{\partial
T}{\partial \xi})
T-\displaystyle\int_{\Omega}T\nabla_{v}T_2-\displaystyle\int_{\Omega}W(v)\frac{\partial
T_2}{\partial
\xi} T\\
&+&\displaystyle\int_{\Omega}\frac{b
P}{p}W(v)T+\displaystyle\int_{\Omega}\frac{ab
P}{p}q_1W(v)T+\displaystyle\int_{\Omega}\frac{ab P}{p}q
W(v_2)T,\qquad\quad\ \quad\quad (6.14)
\end{eqnarray*}
\begin{eqnarray*}
& &\frac{1}{2}\frac{d |q|_2^2}{d
t}+\frac{1}{Rq_1}\displaystyle\int_{\Omega}|\nabla q|^2
+\frac{1}{Rq_2}\displaystyle\int_{\Omega}|q_{\xi}|^2+\frac{\beta_s}{Rq_2}|q|_{\xi=1}|_2^2\\
&=&-\displaystyle\int_{\Omega}(\nabla_{v_1}q+W(v_1)\frac{\partial
q}{\partial \xi})
q-\displaystyle\int_{\Omega}q\nabla_{v}q_2-\displaystyle\int_{\Omega}W(v)\frac{\partial
q_2}{\partial \xi} q.\quad\quad\quad\quad\quad (6.15)
\end{eqnarray*}
Similarly to (6.13), we get
\begin{eqnarray*}
& &\frac{1}{2}\frac{d |T|_2^2}{d
t}+\frac{1}{Rt_1}\displaystyle\int_{\Omega}|\nabla T|^2
+\frac{1}{Rt_2}\displaystyle\int_{\Omega}|T_{\xi}|^2+\frac{\alpha_s}{Rt_2}|T|_{\xi=1}|_2^2\\
&\leq&3\varepsilon\|v\|^2+3\varepsilon\|T\|^2+c|T_2|_4^8(|T|_2^2+|v|_2^2)+c[(|T_{2\xi}|_2^2+1)|\nabla
T_{2\xi}|_2^2+|T_{2\xi}|_2^2]|T|_2^2\\
&+&\displaystyle\int_{\Omega}\frac{b
P}{p}W(v)T+\displaystyle\int_{\Omega}\frac{ab
P}{p}q_1W(v)T+\displaystyle\int_{\Omega}\frac{ab P}{p}qT
W(v_2),\qquad\qquad\quad (6.16)
\end{eqnarray*}
\begin{eqnarray*}
&&\frac{1}{2}\frac{d |q|_2^2}{d
t}+\frac{1}{Rq_1}\displaystyle\int_{\Omega}|\nabla q|^2
+\frac{1}{Rq_2}\displaystyle\int_{\Omega}|q_{\xi}|^2+\frac{\beta_s}{Rq_2}|q|_{\xi=1}|_2^2\\
&\leq&3\varepsilon\|v\|_{H^1(\Omega)}^2+3\varepsilon\|q\|_{H^1(\Omega)}^2
+c[(|q_{2\xi}|_2^2+1)|\nabla
q_{2\xi}|_2^2+|q_{2\xi}|_2^2]|q|_2^2\\
&+&c|q_2|_4^8(|q|_2^2+|v|_2^2).
\quad\quad\quad\quad\quad\quad\quad\quad\quad\quad\quad\quad\quad\quad
\quad\quad\quad\quad\quad\quad\ (6.17)
\end{eqnarray*}
By integration by parts, we have
$$-\displaystyle\int_{\Omega}(\displaystyle\int_{\xi}^1\frac{b
P}{p}\hbox{grad} Td\xi^{'})\cdot
v+\displaystyle\int_{\Omega}\frac{b P}{p}W(v)T=0,\eqno(6.18)$$
$$-\displaystyle\int_{\Omega}(\displaystyle\int_{\xi}^1\frac{ab
P}{p}\hbox{grad} (q_1T)d\xi^{'})\cdot
v+\displaystyle\int_{\Omega}\frac{ab P}{p}q_1W(v)T=0.\eqno(6.19)$$
Similarly to (5.104), we have
\begin{eqnarray*}
|\displaystyle\int_{\Omega}\frac{ab P}{p}q
W(v_2)T|&=&|\displaystyle\int_{\Omega}\displaystyle\int_{\xi}^1\frac{ab
P}{p}\hbox{grad} (q T)d\xi^{'}\cdot v_2|\\
&\leq&c|v_2|_4|q|_4|\nabla T|_2+c|v_2|_4|T|_4|\nabla q|_2\\
&\leq&c|v_2|_4^2(|q|_4^2+|T|_4^2)+\varepsilon|\nabla
T|_2^2+\varepsilon|\nabla q|_2^2\\
&\leq&c|v_2|_4^2(|q|_2^{\frac{1}{2}}\|q\|^{\frac{3}{2}}+|T|_2^{\frac{1}{2}}\|T\|^{\frac{3}{2}})
+\varepsilon|\nabla T|_2^2+\varepsilon|\nabla q|_2^2\\
&\leq&c|v_2|_4^8(|q|_2^2+|T|_2^2)+2\varepsilon\|q\|^2+2\varepsilon\|T\|^2.\quad
\quad\quad\quad(6.20)
\end{eqnarray*}
From (6.13), (6.16)-(6.20), by using (5.12) and (5.13), and
choosing $\varepsilon$ small enough, we obtain
\begin{eqnarray*}
& &\frac{d (|v|_2^2+|T|_2^2+|q|_2^2)}{d
t}+\frac{1}{Re_1}\displaystyle\int_{\Omega}(|\nabla_{e_{\theta}}v|^2+|\nabla_{e_{\varphi}}v|^2+|v|^2)
+\frac{1}{Re_2}\displaystyle\int_{\Omega}|v_{\xi}|^2\\
&+&\frac{1}{Rt_1}\displaystyle\int_{\Omega}|\nabla T|^2
+\frac{1}{Rt_2}\displaystyle\int_{\Omega}|T_{\xi}|^2+\frac{\alpha_s}{Rt_2}|T|_{\xi=1}|_2^2\\
&+&\frac{1}{Rq_1}\displaystyle\int_{\Omega}|\nabla q|^2
+\frac{1}{Rq_2}\displaystyle\int_{\Omega}|q_{\xi}|^2+\frac{\beta_s}{Rq_2}|q|_{\xi=1}|_2^2\\
&\leq&c[|v_2|_4^8+|T_2|_4^8+|q_2|_4^8+|v_{2\xi}|_2^2+(|v_{2\xi}|_2^2+1)\displaystyle\int_{\Omega}(|\nabla_{e_{\theta}}v_{2\xi}|^2
+|\nabla_{e_{\varphi}}v_{2\xi}|^2)]|v|_2^2\\
&+&c(|v_2|_4^8+|T_2|_4^8+|T_{2\xi}|_2^2+(|T_{2\xi}|_2^2+1)|\nabla
T_{2\xi}|_2^2)|T|_2^2\\
&+&c(|v_2|_4^8+|T_2|_4^2+|T_2|_4^4+|q_2|_4^8+|q_{2\xi}|_2^2+(|q_{2\xi}|_2^2+1)|\nabla
q_{2\xi}|_2^2)|q|_2^2. \quad\ (6.21)
\end{eqnarray*}
By Gronwall inequality, Theorem 3.1 and (6.21), we prove Theorem
3.2.

\smallskip
In fact, we have established the following result which is
stronger than Theorem 3.2.

\noindent{\bf Proposition 6.1 (The uniqueness of strong$/$weak
solutions)} \quad Let $U_1$ be a weak solution to the system
$(2.11)-(2.17)$. If there exists a weak solution $U_2$ of the
system (2.11)-(2.17) on the interval $[0,\mathcal{T}]$ with the
same initial conditions, such that
$$U_2\in L^8(0,\mathcal{T};(L^4(\Omega))^4),\ U_{2\xi}\in
L^\infty(0,\mathcal{T};(L^2(\Omega))^4)\cap
L^2(0,\mathcal{T};(H^1(\Omega))^4),$$ Then the solutions $U_1,\
U_2$ coincide on $[0,\mathcal{T}]$.

\section{The existence of universal attractors}
\smallskip
\noindent{\bf Proof of Proposition 3.3} \quad From (5.18), (5.82),
(5.84), (5.98), (5.100), (5.117), (5.118), we know $U\in
L^{\infty}(0,\infty;V)$. By Theorem 3.1 and Theorem 3.2, we can
define the semigroup $\{S(t)\}_{t\geq 0}$ corresponding to the
system (2.11)-(2.16) where $S(t):V\rightarrow V, S(t)U_0=U(t)$. By
(5.16), (5.18), (5.25), (5.32), (5.39), (5.56), (5.64), (5.65),
(5.73), (5.82), (5.83), (5.98), (5.99), (5.117), we prove the
corresponding semigroup $\{S(t)\}_{t\geq 0}$ possesses a bounded
absorbing set $B_\rho$ in $V$, i.e., for any $U_0\in V$, there
exists $t_0$ big enough such that
$$S(t)U_0\in B_\rho, \ \hbox{for any}\ t\geq t_0,$$ where $B_\rho=\{U; \ U\in V,\ \|U\|\leq \rho\}$ and $\rho$ is a
positive constant dependent on $\|Q_1\|_1,\ \|Q_2\|_1$.

In order to prove Theorem 3.4, we need the following property
about the semigroup $\{S(t)\}_{t\geq 0}$.

\smallskip
\noindent{\bf Proposition 7.1} \quad For every $t\geq 0$, the
mapping $S(t)$ is weakly continuous from $V$ to $V$.

\smallskip
\noindent{\bf Proof of Proposition 7.1} \quad Let $\{U_n\}$ be a
sequence in $V$ such that $U_n\rightarrow U$ weakly in $V$.  Then
$\{U_n\}$ is bounded in $V$. By the priori estimates in section 5,
we know that, for every $t\geq 0$,  $\{S(t)U_n\}$ is bounded in
$V$. So we extract a subsequence $\{S(t)U_{n_k}\}$ such that
$S(t)U_{n_k}\rightarrow u$ weakly in $V$. Since the embedding
$V\hookrightarrow L^2(\Omega)\times L^2(\Omega)\times
L^2(\Omega)\times L^2(\Omega)$ is compact, $U_{n_k}\rightarrow U$
strongly in $L^2(\Omega)\times L^2(\Omega)\times L^2(\Omega)\times
L^2(\Omega)$. By (6.21), we obtain that $S(t)U_{n_k}\rightarrow
S(t)U$ strongly in $L^2(\Omega)\times L^2(\Omega)\times
L^2(\Omega)\times L^2(\Omega)$. Then $u=S(t)U$. Therefore, the the
sequence $\{S(t)U_n\}$ satisfies: $S(t)U_{n_k}\rightarrow S(t)U$
weakly in $V$. Proposition 7.1 is proved.

\smallskip
\noindent{\bf Proof of Theorem 3.4} \quad With Proposition 3.3 and
Proposition 7.1, we know that the proof of Theorem 3.4 is similar
to that of Theorem I.1.1 in \cite{T1}. We only need replace
"$\rightarrow \ \hbox{strongly in}\  H$" by "$\rightarrow \
\hbox{weakly in} \ V$" in the course of proof of Theorem I.1.1 in
\cite{T1}. So the detail of proof for Theorem 3.4 is omitted here.

\bigskip
\noindent {\bf Acknowledgments.} The second author would like to
express my heartful thanks to Prof. Songying  Li for the useful
comments and suggestions. The work was supported in part by the
NNSF of China grants No 90511009.

\bigskip

\end{document}